\documentclass[journal,twoside]{IEEEtran}

\usepackage[pdftex]{graphicx}
\graphicspath{{../pdf/}{../jpeg/}}
\DeclareGraphicsExtensions{.pdf,.jpeg,.png}
\usepackage{amssymb,amsmath,url}
\usepackage{bm} 
\usepackage{multirow}
\usepackage[dvipsnames]{xcolor}
\usepackage{booktabs}	
\usepackage{epstopdf}
\usepackage{subfigure}
\usepackage{bigfoot}		
\usepackage{algpseudocode}
\usepackage{pifont}
\usepackage{tikz}
\usepackage{comment}
\usepackage{cite}
\usepackage{lipsum}
\usepackage{algorithm} 
\usepackage[pdftex]{graphicx}
\usepackage{array}
\usepackage{enumitem}

\algrenewcommand\algorithmicrequire{\textbf{Input:}}
\algrenewcommand\algorithmicensure{\textbf{Output:}}
\algrenewcommand\algorithmicforall{\textbf{For}}

\newtheorem{theorem}{Theorem}
\newtheorem{lemma}{Lemma}
\newtheorem{proposition}{Proposition}
\newtheorem{corollary}{Corollary}

\allowdisplaybreaks
\newtheorem{definition}{Definition}

\newtheorem{remark}{Remark}

\newtheorem{example}{Example}

\newcommand{\R}{\mathbb{R}}

\newcommand*{\QEDA}{\hfill\ensuremath{\blacksquare}}

\title{
	Actuator Placement under Structural Controllability using Forward and Reverse Greedy Algorithms
}

\author{Baiwei Guo$^\ast$\thanks{$^\ast$These authors contributed equally to this work.\\ \indent The work of Karaca and Kamgarpour was gratefully funded by the European Union ERC Starting Grant CONENE. The work of Summers was sponsored by the Army Research Office and was accomplished under Grant Number: W911NF-17-1-0058.}, Orcun Karaca$^\ast$, Tyler Summers,  Maryam Kamgarpour\thanks{B. Guo is with the Automatic Control Laboratory, EPFL, Switzerland. {\tt baiwei.guo@epfl.ch}\ O. Karaca, M. Kamgarpour are with the Automatic Control Laboratory, D-ITET, ETH Z\"{u}rich, Switzerland. {\tt \{okaraca, mkamgar\}@ethz.ch}\ T. Summers is with the Dept. of Mech. Eng., UT Dallas, Richardson, TX, USA. {\tt tyler.summers@utdallas.edu}.}
}%
\allowdisplaybreaks
\newcommand\scalemath[2]{\scalebox{#1}{\mbox{\ensuremath{\displaystyle #2}}}}
\begin{document}

\maketitle 
\begin{abstract}\noindent Actuator placement is an active field of research which has received significant attention for its applications in complex dynamical networks. In this paper, we study the problem of finding a set of actuator placements minimizing the metric that measures the average energy consumed for state transfer by the controller, while satisfying a structural controllability requirement and a cardinality constraint on the number of actuators allowed. 
As no computationally efficient methods are known to solve such combinatorial set function optimization problems, two greedy algorithms, forward and reverse, are proposed to obtain approximate solutions. We first show that the constraint sets these algorithms explore can be characterized by matroids. 
We then obtain performance guarantees for the forward and reverse greedy algorithms applied to the general class of matroid optimization problems by exploiting properties of the objective function such as the submodularity ratio and the curvature. Finally, we propose feasibility check methods for both algorithms based on maximum flow problems on certain auxiliary graphs originating from the network graph. Our results are verified with case studies over large networks.
\end{abstract}

\section{Introduction}

\IEEEPARstart{M}{any} large-scale complex dynamical networks, such as those arising in power grids \cite{summers2015submodularity}, biological networks \cite{muller2011few} and industrial systems \cite{banerjee1995control} necessitate a resilient and efficient operation under dynamic and uncertain environments. Hence, there has been a surge of interest to study controller design in such large-scale networks~\cite{van2001review,liu2011controllability,pasqualetti2014controllability,olshevsky2014minimal,pequito2015framework,Ali2017MinimalReachablity,clark2017submodularity,romero2018actuator,summers2018performance,pequito2017robust,tzoumas2019robust}. A fundamental design problem is that of actuator placement in which the goal is to select a subset from a finite set of possible placements for actuators to optimize a desired network performance metric.

Variants of the actuator placement problem have been shown to be NP-hard in general, see~\cite{olshevsky2014minimal,EffortBounds, Ali2017MinimalReachablity}. Thus, it is desirable to obtain scalable algorithms with provable suboptimality bounds. Earlier studies have adopted the forward greedy algorithm. This algorithm extends the actuator set with the most beneficial actuator iteratively to derive an approximate solution~\cite{summers2015submodularity}. Under a submodular network performance metric and a cardinality constraint on the number of actuators, the forward greedy algorithm is shown to enjoy a provable performance guarantee~\cite{nemhauser1978analysis}. However, some metrics do not exhibit submodularity including the metric in this work, that is, the average energy required to reach any arbitrary direction of the state space~\cite{summers2017performance}. To alleviate this issue, submodularity has been extended to weak submodularity using the notion of submodularity ratio, quantifying how close a function is to being submodular~\cite{Das:2011:SMS:3104482.3104615,bian2017guarantees}. Given this ratio, it is possible to derive a performance guarantee for the forward greedy algorithm applied to a larger class of performance~metrics~\cite{summers2017performance}.

Nonetheless, the guarantees above are restricted to problems subject to simple cardinality constraints. Given a cardinality constraint, the resulting actuator {set} might not be capable of moving the system {over} the entire state space, that is, might not render the system controllable. 
To address this issue, we need to include controllability as a constraint. {However, to the best of our knowledge, there is no approach to quantify the forward greedy algorithm's performance with a nonsubmodular metric and a controllability constraint, nor to ensure feasibility of the iterates of the greedy algorithm in such problems.} On the other hand, structural controllability constraints have been well-studied. This controllability concept exploits only the graphical interconnection structure of the dynamical system \cite{SturcturalcontrollabilityTai,shields1976structural,dion2003generic,liu2011controllability,ramos2020structural}. Structurally controllable systems are those controllable after a slight perturbation of the system parameters corresponding to the fixed set of edges in the underlying network graph. The authors in~\cite{clark2012leader} have studied a leader selection problem to obtain a structurally controllable system while minimizing a submodular objective function. The structural controllability constraint arising in the leader selection problem is proven to give rise to a matroid constraint enabling the application of the forward greedy algorithm~\cite{clark2012leader}. However, the leader selection problem is different from the actuator placement problem. The former selects a set of leader nodes whose states can arbitrarily be dictated to steer {the remaining nodes to desired states},
while the latter does not permit the states to be dictated arbitrarily; instead, it selects a set of actuators which can influence all of the states through the dynamics. {Hence, this paper pays special attention to formulating the structural controllability constraints of the actuator placement problem as a matroid constraint by proving the equivalence of this concept in both the leader selection problems and the actuator placement problems.}

Given a matroid, \cite{fisher1978analysis} derives a performance guarantee for the forward greedy algorithm when optimizing a submodular objective. However, past work has not successfully derived performance guarantees for optimizing weakly submodular objective functions, such as the aforementioned average energy consumption metric, subject to a matroid. The first goal is to obtain a guarantee for this setting. In Appendix~\ref{app:tableg}, we discuss relevant existing guarantees from~\cite{nemhauser1978analysis,fisher1978analysis,bian2017guarantees,conforti1984submodular,pmlr-v80-chen18b,chamon2019matroid,sviridenko2017optimal}.
\looseness=-1

An inherent drawback of the forward greedy algorithm is that any performance guarantee has to involve the objective function evaluated at the empty set as the reference value, since the actuator set expands starting from the empty set. This reference value is in general large for the average energy consumption metric, or even infinite~\cite{summers2017performance}, and it plays a great role towards the tightness of the guarantee. In addition, many works have reported the lack of ability of the forward greedy to correct errors made in earlier steps~\cite{zhang2011adaptive,tropp2004greed}.
An alternative is to adopt the reverse greedy, which excludes the least beneficial actuator iteratively starting from the full set. In this case, any potential performance guarantee would instead involve the objective function evaluated at the full set, which is in general small for the performance metric considered in this work. 
\looseness=-1

Among the applications of the reverse greedy algorithm, \cite{chrobak2006reverse} studied the special setting of metric $k$-median problem and this algorithm is {shown} to have a better performance than the forward greedy algorithm. The work of~\cite{il2001approximation} provides a guarantee for minimizing a supermodular decreasing function under cardinality constraints by exploiting a notion of function steepness, while \cite{il2006performance} extends this analysis to account for comatroid constraints.\footnote{Comatroid is the complementary notion of a matroid, see \cite{il2006performance,il2003hereditary}} Our paper in~\cite{karaca2019} provides a counterexample to the performance guarantee obtained in~\cite{il2006performance}, and explains where the mistake originates from in their proof. Nevertheless, none of the problem settings can generalize the problem of actuator placement considered in this work.
This is because, in addition to involving matroid constraints, via a reformulation, the objective function of our problem will be shown to exhibit weak supermodularity, which will be characterized by the notion of curvature~\cite{sviridenko2017optimal,bian2017guarantees}.
To the best of our knowledge, there is no performance guarantee for the reverse greedy algorithm applicable to optimizing weakly submodular and weakly supermodular objective functions (defined by submodularity ratio and curvature, respectively) subject to matroid constraints.
\looseness=-1

Our main contributions are as follows.

(i) We show that the minimization of the average energy consumption metric under structural controllability constraints can be reformulated as the maximization of a strictly increasing weakly submodular function subject to matroid constraints, see Lemma~\ref{lmm: monotonicity}, Proposition~\ref{prop: matroid}, and Problem~\eqref{eq:mainproblem}. 

(ii) We obtain a performance guarantee for the forward greedy algorithm applied to this general class of matroid optimization problems, see Theorem~\ref{thm:upperlimitFormatroidconstraints}.\footnote{{Theorems~\ref{thm:upperlimitFormatroidconstraints} and \ref{thm:reverseGuarantee} could be of independent interest for researchers working on greedy algorithms.} Preliminary results concerning the forward greedy---(i) and (ii) above---were presented in a conference paper in~\cite{guo2019actuator}. This paper significantly extends that work by contributions (iii) to (vii), and utilizes the newly introduced greedy notions of the curvature and the submodularity ratio. }
    
(iii) We show that the actuator placement problem has another reformulation as the minimization of a strictly increasing, weakly submodular, and weakly supermodular function subject to matroid constraints and a cardinality lower bound, see Lemma~\ref{lem:weakrev}, Proposition~\ref{prop: matroid2}, and Problem~\eqref{eq:mainproblemanotherform}. This reformulation allows us to implement the reverse greedy algorithm.
    
(iv) 
For the reverse greedy algorithm, we obtain a performance guarantee employing both notions, see Theorem \ref{thm:reverseGuarantee}. 
    
(v) The average energy consumption metric is well-defined only if we introduce a metric-modifying parameter~\cite{EffortBounds}.
To this end, we design an algorithm with a provable performance to pick such parameters, see Proposition~\ref{prop:epsilonalgproof} and Algorithm~\ref{alg:epsilon}.

(vi) For both algorithms, we show that the matroid feasibility checks for the actuator placement
    can be done efficiently by translating them into maximum flow problems over certain auxiliary graphs, see Propositions~\ref{prop: feasibility},~\ref{prop: feasibilityreverse}. These results extend \cite{liu2011controllability} which associates structural controllability with the existence of a perfect matching. We also provide a counterexample to a feasibility check in~\cite{clark2012leader} for the leader selection problem.
    \looseness=-1
    
{Finally, we provide numerical case studies with models based on randomly generated networks and a large power grid.} As an additional insight, we demonstrate that the forward greedy algorithm tends to pick higher degree actuators when compared to the optimal and the reverse greedy solutions.
\looseness=-1
 
In the remainder, Section~\ref{sec:mech} introduces the problem formulation and preliminaries. Sections~\ref{sec:Forwardgreedyalgorithm} and~\ref{sec:reverseGreedy} apply the forward and the reverse greedy algorithms, respectively, and obtain guarantees. Section~\ref{sec:FeasibilityCheck} proposes a method to pick a metric-modifying parameter and feasibility check methods for greedy algorithms. Numerical studies are presented in Section~\ref{sec:num}. 
\looseness=-1

\section{Problem Formulation and Preliminaries}\label{sec:mech}

\subsection{Problem formulation}\label{sec:2a}
Consider a linear system with state vector $x\in \R^n$. To each state variable $x_i\in \R$, we associate a node $v_i\in V:= \{v_1, \ldots, v_n\}$. A control input $u_i\in \R$ can be exerted at each node~$v_i$. Given a set $S\subset V$ chosen as the actuator set, the system dynamics can be written as
\begin{equation}
\begin{aligned}
\label{eq: systemmodel}
\dot{x} = Ax + B(S)u.
\end{aligned}
\end{equation}
Above, $B(S)= \text{diag}(\bm {1}(S))\in \R^{n\times n}$, where $\bm {1}(S)$ denotes a vector of size $n$ whose $i$th entry is $1$ if $v_i$ belongs to $S$ and $0$ otherwise. {Weights of the entries in $B$ are not decision variables, as opposed to the works of~\cite{siami2018deterministic}.}  Let $G=(V,E)$ denote a directed graph relating to system~\eqref{eq: systemmodel} with nodes $V$ and edges $E$, where the edge $(v_j,v_i)\in E$ if $(A)_{ij}\neq 0$. Similar to several previous studies on structural controllability, e.g.,  \cite{clark2012leader,clark2017submodularity}, throughout the paper, we assume that {$G$ is strongly connected, which will be discussed in Section~\ref{sec:feasibilitycheck}.}

{The pair $(A,B(S))$ is called \textit{controllable} if
for all $x_0,x_1\in\R^n$ and $T>0$ there exists a control input $u:[0,\, T]\rightarrow\R^n$ that steers the system from $x_0$ at $t=0$ to $x_1$ at $t=T$.} For linear time-invariant systems, controllability can be verified by the rank of the controllability matrix $P=\begin{bmatrix}
   B(S) & AB(S) & \cdots & A^{n-1}B(S)
\end{bmatrix}\in \R^{n\times n^2}$. However, the entries in $A$ are generally not exactly~known
but only approximately determined with small errors using system identification techniques. Moreover, when dealing with large-scale networked systems, it is often the case that we can only rely on the topology but not on the particular weights~\cite{romero2018actuator}. Motivated by these particularities, we consider structural controllability.

\begin{definition}
\label{def: structuralcontrollability}
$(A,B)$ and $(\hat{A},\hat{B})$ with $A,B,\hat{A},\hat{B}$ $\in \R^{n\times n}$ are said to have the same structure if matrices $[A\text{ }B]$ and $[\hat{A} \text{ }\hat{B}]$ have zeros at the same entries. Given $S\subset V$, $(A,B(S))$ is \textit{structurally controllable} if there exists a controllable pair $(\hat{A},\hat{B})$ having the same structure as $(A,B(S))$.
\end{definition}

As it turns out, structural controllability is a generic property, that is, the pair $(A,B(S))$ is structurally controllable if and only if almost all of the pairs with the same structure are controllable~\cite{dion2003generic}. This implies that whenever $(A,B(S))$ is not controllable but structurally controllable, it is possible to slightly perturb the entries to ensure controllability~\cite{SturcturalcontrollabilityTai}.
Observe that structural controllability depends on the positions of the nonzero entries. Later, this will allow us to determine this property by the graph $G$ relating to the system.
\looseness=-1

Even if a system is controllable, an unacceptably large amount of energy might be needed to reach a desired state. 
Specifically, the work in~\cite{pasqualetti2014controllability} shows that if the number of actuators is kept constant, then certain controllable systems are practically uncontrollable since the energy consumption grows at least exponentially with the number of states~$n$. Hence, it is crucial to minimize this energy consumption. The minimum energy required to steer the system from zero at $t=0$ to $x\in\R^n$ at $t=T$ is given by $x^\top W^{-1}_T(S)x,$ where $W_T(S)=\int_{0}^{T} e^{A\tau}B(S) B^\top (S) e^{A^\top\tau} d\tau$ is the controllability Gramian. To obtain an expression independent of the initial state $x$, calculate the average energy required over the unit sphere, $||x||_2=1$, as $F(S):=\text{tr}(W_T^{-1}(S))$. This expression is well-defined only when the set $S$ renders the system controllable. Inspired by~\cite{EffortBounds}, we introduce a small positive number $\epsilon\in \R_+$ to handle uncontrollable actuator sets and propose the metric $F_\epsilon:2^V\rightarrow\R_+,$
\begin{equation}
\label{eq:F(S)}
F_\epsilon(S) =  \text{tr}((W_T(S)+\epsilon I)^{-1}), \text{ }\forall S\subset V.
\end{equation}
In Section~\ref{sec: algorithmforproperep}, we discuss the choice of $\epsilon$.

To make a system easier to control, we seek a set $S\subset V$ minimizing the metric above. Since in a large-scale network, the number of actuators allowed is in general limited, we consider a cardinality bound of $K\in \mathbb{N}$ on the actuators. Additionally, we require that the actuators render the system structurally controllable. Our main problem is formulated as
\begin{equation}
\begin{aligned}
& \min_{S\subset V}
& & F_\epsilon(S)
\\ 
&\ \mathrm{s.t.} & & |S|\leq K,\, \text{$(A,B(S))$ is structurally controllable.}
\end{aligned}
\label{eq:mainproblem2}
\end{equation}
Assume $K$ is large enough to ensure feasibility. {In Section \ref{sec:feasibilitycheck}, we discuss how to determine the smallest such $K$.}
Problem \eqref{eq:mainproblem2} is a combinatorial optimization, and to the best of our knowledge, no computationally feasible solution method has ever been proposed. {Existing works have either studied additive/modular objectives~\cite{pequito2013structured,pequito2015framework,romero2018actuator,doostmohammadian2018structural} (e.g., actuator installation costs, minimizing $K$ directly), or included only cardinality constraints~\cite{summers2015submodularity,summers2017performance}.} Notice that neither our objective is additive nor we have only cardinality constraints. Later, we will adopt efficient heuristics to derive solutions.
\looseness=-1

\subsection{Preliminaries}

We first introduce widely adopted notions for the properties of set functions and set constraints.

\subsubsection{Properties of set functions}

Given a ground set $V$ and a set function $f:2^V\rightarrow \R$, we say $f$ is \textit{(strictly) increasing} if $f(S_1)\leq$($<$)$ f(S_2)$ for any $S_1\subsetneqq S_2\subset V$. If $-f$ is (strictly) increasing, we say $f$ is \textit{(strictly) decreasing}. 
For an increasing set function, the marginal gain from the addition of a certain element $v\in V$ to a set $S\subset V$ varies for different $S$. {For many set functions in practical problems the marginal gain diminishes as $S$ expands, see the examples in \cite{Krause2005Near,bach2013learning}.} Submodularity describes this property and submodularity ratio describes how far a nonsubmodular function is from being submodular. For the following, denote the marginal gains by $\rho_{U}(S):=f(S\cup U)-f(S),\ \forall S,U\subset V.$ For notational simplicity, we use $v$ and $\{v\}$ interchangeably for singleton sets.
\looseness=-1

\begin{definition}\label{def:submodularityratio}
For an increasing function $f \colon 2^V \to \R$, the \textit{submodularity ratio} is the largest $\gamma \in \R_+$ such that
$\gamma\rho_v(S\cup U)\leq  \rho_v(S),\text{ } \forall S,U \subset V, \text{ }\forall v \in V\setminus(S\cup U).$
It can be verified that $\gamma \in [0,1].$
A set function $f$ with submodularity ratio~$\gamma$ is called $\gamma$-submodular. A $\gamma$-submodular set function is said to be \textit{submodular} if $\gamma=1$ and \textit{weakly submodular} if $0<\gamma<1$.
\end{definition}

In Appendix~\ref{apd: definitions of submodularity ratio}, we connect Definition~\ref{def:submodularityratio} with another existing notion of submodularity ratio and discuss the necessity of introducing this notion as per Definition~\ref{def:submodularityratio} for the guarantee derived for the forward greedy algorithm in Section~\ref{sec:Forwardgreedyalgorithm}.

Other than submodularity, another widely-used notion is supermodularity, that is, the marginal gain from the addition of $v\notin S$ to the set $S$ increases as $S$ expands. By introducing supermodularity and the curvature, that is, how far a nonsupermodular function is from being supermodular, we obtain a more precise description on how the marginal gains change. 
\begin{definition}
\label{def:curvature}
For an increasing function $f:2^V \to \R$, the \textit{curvature} is the smallest $\alpha \in \mathbb{R}_+$ such that
$\rho_v(S\cup U)\geq  (1-\alpha)\rho_v(S),\text{ } \forall S,U \subset V, \text{ }\forall v \in V\setminus(S\cup U).$ It can be verified that $\alpha \in [0,1].$
Function $f$ with curvature $\alpha$ is called $\alpha$-supermodular. 
 An $\alpha$-supermodular function is \textit{supermodular} if $\alpha=0$ and \textit{weakly supermodular} if $0<\alpha<1$.
\end{definition}

To see how submodularity ratio and curvature are related, notice that for an increasing set function $f$ the submodularity ratio $\gamma$ and the curvature $\alpha$ satisfy
{\medmuskip=1.25mu\thinmuskip=1.25mu\thickmuskip=1.25mu	\begin{equation}
\label{eq:extentionofGammaAlphaDef}
\gamma = \min_{\small\substack{S,U,\\v\in V\setminus(S\cup U)}}\frac{\rho_v(S)}{\rho_v(S\cup U)}\leq \max_{\small\substack{S,U,\\v\in V\setminus(S\cup U)}}\frac{\rho_v(S)}{\rho_v(S\cup U)}=\frac{1}{1-\alpha}.    
\end{equation}}

\subsubsection{Properties of set constraints}

{Many combinatorial optimization problems from the literature are subject to constraints that are more complex than simple cardinality constraints, see the examples in~\cite{krause2014submodular,tzoumas2018resilient}. Among those, we introduce matroids since they will generalize reformulations of the constraints found in Problem \eqref{eq:mainproblem2}, and they allow performance guarantees for greedy algorithms~\cite{edmonds1971matroids}.}
\begin{definition}
\label{def:matroid}
A \textit{matroid} $\mathcal{M}$ is an ordered pair $(V,\mathcal{F})$ consisting of a ground set $V$ and a collection $\mathcal{F}$ of subsets of $V$ which satisfies (i) $\emptyset \in \mathcal{F}$, (ii) if $S,S' \in \mathcal{F}$ and $S'\subset S$, then $S'\in \mathcal{F}$, (iii) if $S_1$,$S_2\in \mathcal{F}$ and $|S_1|<|S_2|$, there exists $v \in S_2\setminus S_1$ such that $v\cup S_1\in \mathcal{F}$. Every set in $\mathcal{F}$ is called \textit{independent}, and maximum independent sets refer to those with the largest cardinality.
\end{definition}

To adopt the reverse greedy algorithm, an additional concept will be required, that is, the dual of a matroid.
\begin{definition}
\label{def:dual of a matroid}
Given a matroid $(V,\mathcal{F})$, let $\mathcal{F}^* = \{U\text{ }|\text{ }\exists \text{ a maximum independent set}$ $ M\in \text{  $\mathcal{F}$ such}$ $\text{that } U\subset V\setminus M\}$. The pair $(V,\mathcal{F}^*)$ is the \textit{dual} of the matroid $(V,\mathcal{F})$.
\end{definition}

We characterize its structure in the following lemma. 

\begin{lemma}
\label{prop:dualmatroid}
{The pair $(V,\mathcal{F}^*)$, the dual of a matroid $(V,\mathcal{F})$, is also a matroid.}
\end{lemma}
\begin{IEEEproof}
Suppose $\{M_i\}_{i=1}^q$ is the collection of all maximum independent sets in matroid $(V,\mathcal{F})$. From \cite[Ch. 2]{welsh2010matroid} we have that $\{V\setminus{M_i}\}_{i=1}^q$ defines a collection of all maximum independent sets for another matroid denoted by $(V,\tilde{\mathcal{F}})$. In the following, we prove that $\mathcal{F}^* = \tilde{\mathcal{F}}$.
For any $U\in \mathcal{F}^*$, there exists $M$, a maximum independent set in $\mathcal{F}$, such that $U\subset V\setminus M$. Since $V\setminus M\in \tilde{\mathcal{F}}$ and $(V,\mathcal{\tilde{F}})$ is a matroid, the set $U$ also belongs to $\tilde{\mathcal{F}}$ from property (ii) in Definition~\ref{def:matroid}. Conversely, if $U\in \tilde{\mathcal{F}}$, according to property {(iii)} in Definition~\ref{def:matroid}, $U$ is a subset of some maximum independent set in $\tilde{\mathcal{F}}$. Consequently, there exists a maximum independent set $M\in \mathcal{F}$ such that $U\subset V\setminus{M}.$ Thus, $U\in \mathcal{F}^*.$ This concludes that $\mathcal{F}^*=\mathcal{\tilde{F}}$ and thus $(V,\mathcal{F}^*)$ is also a matroid. 
\end{IEEEproof} 

{\cite[Ch. 2]{welsh2010matroid} defines the dual concept as $(V,\tilde{\mathcal{F}})$, and the proof above verifies that $(V,\mathcal{F}^*)$ we have in Definition~\ref{def:dual of a matroid} is an equivalent reformulation. This reformulation will help us present the proof of Proposition~\ref{prop: matroid2} in a more clear way.}

\section{Forward Greedy Algorithm}
\label{sec:Forwardgreedyalgorithm}

In the following, we reformulate Problem \eqref{eq:mainproblem2} as the maximization of a strictly increasing weakly submodular function subject to matroid constraints. We then obtain a guarantee for a forward greedy algorithm over matroid constraints.
\looseness=-1

\subsection{Properties of the objective}
Intuitively, with more input nodes, system~\eqref{eq: systemmodel} would be easier to control and thus the metric $F_{\epsilon}$ in~(\ref{eq:F(S)}) would be smaller. This intuition can be readily verified as follows.

\begin{lemma}
\label{lmm: monotonicity}
The metric $F_\epsilon=\text{tr}((W_T(S)+\epsilon I)^{-1})$ satisfies the following statements: (i) $F_\epsilon$ is strictly decreasing, (ii) $-F_\epsilon$ is weakly submodular with submodularity ratio $\gamma_\epsilon^\mathsf{f}$.
\end{lemma}

The proof is relegated to Appendix \ref{apd: proofofmonotonicity}. Together with the fact that the structural controllability is preserved under actuator set expansion, Lemma~\ref{lmm: monotonicity} implies that the optimal solution to Problem~\eqref{eq:mainproblem2} should contain exactly $K$ nodes.

\subsection{Reformulation of the constraint set}
In combinatorial optimization problems with only cardinality constraints, the forward greedy algorithm starts from the empty set and at $t^{\text{th}}$ iteration, adds the most marginally beneficial node $v_t^\mathsf{f}$ to the actuator set. It terminates when the cardinality of the actuator set is~$K$. When applied to Problem~\eqref{eq:mainproblem2}, this method might return an actuator set under which the system is not structurally controllable. To this end, we need to restrict the greedy iterates $S^t=\{v_1^\mathsf{f},\ldots, v_t^\mathsf{f}\}$, for $t=1,\ldots, K$, such that the set $S^K$ returned by the forward greedy algorithm is guaranteed to satisfy structural controllability.
\looseness=-1

Since the optimal solution to Problem~\eqref{eq:mainproblem2} contains exactly $K$ nodes, we define $\mathcal{C}_K=\{S\subset V\,|\,|S| = K$ and the system is structurally controllable under $S\}$ and rewrite Problem \eqref{eq:mainproblem2} as the minimization of $F_\epsilon$ over the set collection $\mathcal{C}_K$. In the procedure of the forward greedy algorithm, the set $S^t$ has to be a subset of some set in $\mathcal{C}_K$, since otherwise the greedy solution $S^K$ would not belong to $\mathcal{C}_K$. Thus, define $\mathcal{\tilde{C}}_K=\{\Omega\,|\,\exists S \in \mathcal{C}_K \text{ such that } \Omega \subset S \}$ and reformulate~\eqref{eq:mainproblem2}~as   
\begin{equation}
 \max_{S\subset V}\
 -F_\epsilon(S)
\ \ \mathrm{s.t.}\ S\in \mathcal{\tilde{C}}_K.
\label{eq:mainproblem}
\end{equation}
The strict monotonicity of $-F_\epsilon$ ensures that the optimal solution to Problem~\eqref{eq:mainproblem} coincides with that of Problem~\eqref{eq:mainproblem2}. As such, we consider solving Problem \eqref{eq:mainproblem} as an equivalent characterization of Problem~\eqref{eq:mainproblem2}. 

{Next, we show that the feasible region of Problem~\eqref{eq:mainproblem}  characterizes a matroid, which will allow us to derive performance guarantees for the greedy solution $S^K$.}
\begin{proposition}
\label{prop: matroid}
$\mathcal{M}=(V,\tilde{\mathcal{C}}_K)$ is a matroid. 
\end{proposition}

{To prove this, we establish the equivalence between structural controllability of $(A,B(S))$ in {Problem \eqref{eq:mainproblem2}} and structural controllability of the system with the set $S$ chosen as a leader set in a corresponding leader selection problem. We then invoke a result from~\cite{clark2012leader} proving the matroid structure of the structural controllability constraints in leader selection problems. For the details and additional discussions, we kindly refer to the proof in Appendix \ref{apd: matroid}. This equivalence result will later be utilized in Sections~\ref{sec:feasibilitycheck} and~\ref{sec:reversefeasibilitycheck} to bring in results from the leader selection literature.}
\looseness=-1

We now restrict the iterates of the forward greedy algorithm to lie in the set collection $\tilde{\mathcal{C}}_K$. 
As a remark, given $S^t\in \tilde{\mathcal{C}}_K$ for $t<K$, by Definition~\ref{def:matroid}, we can always find a node $v\in V\setminus S^t$ such that $S^t\cup v\in  \tilde{\mathcal{C}}_K$, as long as $\mathcal{C}_K\neq\emptyset.$ Therefore, it is guaranteed that at iteration~$K$ we obtain an actuator set in $\mathcal{C}_K$.

\subsection{Performance guarantee}
\label{sec: performance guarantee}
In the previous section, we showed that the objective function $-F_\epsilon$ is  $\gamma_\epsilon^\mathsf{f}$-submodular in Lemma~\ref{lmm: monotonicity} and the feasible region $\tilde{\mathcal{C}}_K$ characterizes a matroid in Proposition~\ref{prop: matroid}. Thus, Problem~\eqref{eq:mainproblem} falls into the following class of optimization problems:
\begin{equation}
\begin{aligned}
& \max_{S\subset V}
& & f(S), \text{strictly increasing and $\gamma$-submodular}
\\ 
& \ \mathrm{s.t.} & & S\in \mathcal{F}, \text{ where $\mathcal{M}=(V,\mathcal{F})$ is a matroid},
\end{aligned}
\label{eq:matroidconstrained optimization}
\end{equation}
where the cardinality of any maximum independent set in $\mathcal{F}$ is~$K$. Let $S^*$ denote its optimal solution.

The forward greedy over a matroid was first introduced in \cite{fisher1978analysis} for submodular objectives. 
This algorithm is presented in Algorithm \ref{alg:ALG4}. At the $t^\text{th}$ iteration, we check the feasibility of the node with the largest marginal gain in $V\setminus S^{t-1}$. If the actuator set obtained by adding this node to $S^{t-1}$ does not belong to $\mathcal{F}$, we exclude the node from consideration. Among the remaining ones, we check the feasibility of the node with the largest marginal gain until a feasible node $v^\mathsf{f}_t$ is found. Then $S^t=\{v^\mathsf{f}_t\}\cup S^{t-1}$ is the actuator set returned by the $t^\text{th}$ iteration. The final actuator set is $S^\mathsf{f}:=S^K$. The feasibility check ensures that $S^t\in \mathcal{F}$ and hence $S^\mathsf{f}$ belongs to~$\mathcal{F}$.

We use $U^t\subset V$ for $0\leq t\leq K-1$ to denote all the nodes having been considered by the feasibility check before $v^\mathsf{f}_{t+1}$. We define the marginal gains as $\rho_t = f(S^{t+1})-f(S^t)$. 

\begin{algorithm}[t]
        \caption{Forward Greedy Algorithm over Matroid}
        \label{alg:ALG4}
        {\begin{algorithmic}
        \Require set function $f$, ground set $V$ and matroid $(V,\mathcal{F})$
        \Ensure actuator set $S^\mathsf{f}$
        \Function {ForwardOverMatroid}{$f, V,\mathcal{F}$}
        \State {$S^0 = \emptyset$, $U^0= \emptyset$, $t=1$}
        \While{$U^{t-1}\neq V$ { and } $|S^{t-1}|<K$}
            \State ${i^*(t)}= \arg\max_{i\in V\setminus U^{t-1}}\rho_i(S^{t-1})$
            \If {$S^{t-1}\cup\{i^*(t)\}\notin\mathcal{F}$}
        \State $U^{t-1}\gets U^{t-1}\cup\{i^*(t)\}$
        \Else
        \State $\rho_{t-1} \gets\rho_{i^*(t)}(S^{t-1})$ and $v^\mathsf{f}_t = i^*(t)$
        \State $S^{t} \gets S^{t-1}\cup\{v^\mathsf{f}_t\}$ and $U^{t}\gets U^{t-1}\cup\{v^\mathsf{f}_t\}$
        \State $t\gets t+1$
    \EndIf    
            \EndWhile
            \State $S^\mathsf{f} \gets S^{t-1}$
        \EndFunction
        \end{algorithmic}}
    \end{algorithm}

Using the matroid structure and the submodularity ratio, we can state our first main result as follows.
\begin{theorem}
If Algorithm \ref{alg:ALG4} is applied to Problem (\ref{eq:matroidconstrained optimization}), then 
\begin{equation}
\begin{aligned}
\frac{f(S^\mathsf{f})-f(\emptyset)}{f(S^*)-f(\emptyset)}\geq\frac{\gamma^3}{\gamma^3+1}.
\end{aligned}
\label{eq: performanceguaranteeMatroidOpt}
\end{equation}
\label{thm:upperlimitFormatroidconstraints}
\end{theorem}

The proof is relegated to Appendix~\ref{apd: theorem1}. The idea of the proof extends the work in \cite{fisher1978analysis}, which derives a performance guarantee for matroid optimization featuring a submodular objective. When $\gamma = 1$, the guarantee in \eqref{eq: performanceguaranteeMatroidOpt} coincides with that of \cite{fisher1978analysis}, derived for a submodular $f$. 
As a remark, for Problem~\eqref{eq:matroidconstrained optimization}, another performance guarantee is offered by \cite{pmlr-v80-chen18b} {but in expectation for a randomized algorithm}. We refer to Appendix \ref{apd: definitions of submodularity ratio} for a comparison of these two guarantees.\footnote{The works in \cite{conforti1984submodular} and \cite{bian2017guarantees} utilize also the curvature to derive performance guarantees for the forward greedy applied to cardinality constrained problems. Exploiting this notion for matroid constraints is part of our ongoing work. }

Given any function $f$, it is difficult to derive its submodularity ratio because the computation in Definition~\ref{def:submodularityratio} involves {${\Omega}(2^n)$} inequalities.
In the proof of Theorem \ref{thm:upperlimitFormatroidconstraints}, only a subset of these inequalities are utilized. Via this observation, the following corollary proposes a computationally more efficient approach.
\looseness=-1

\begin{corollary}\label{cor:1}
Let $\gamma^\mathsf{fg}$ be the largest {$\hat\gamma$} that satisfies
(a) $f(S\cup S^\mathsf{f})-f(S^\mathsf{f})\leq {\hat{\gamma}^{-1}}\sum_{j\in S\setminus S^\mathsf{f}}\rho_j(S^\mathsf{f})$ for any $S$ with $|S|=K$, (b) $\rho_j(S^\mathsf{f})  \leq {\hat{\gamma}^{-1}}\rho_j(S^{t-1}), \forall t\leq K\text{, } \forall j\in V$,
(c) $ f(S^{i_2+1})-f(S^{i_2})
	\leq {\hat{\gamma}^{-1}}(f(S^{i_1}\cup\{v^\mathsf{f}_{i_2+1}\})-f(S^{i_1})), $ for any $i_1<i_2$. Then, $\gamma^{\mathsf{fg}}$ is called the greedy submodularity ratio for the forward greedy algorithm, with $\gamma^{\mathsf{fg}}\geq \gamma$, and 
\begin{equation}
\label{eq: greedyversionPerformanceGuarantee}
\frac{f(S^\mathsf{f})-f(\emptyset)}{f(S^*)-f(\emptyset)}\geq\frac{(\gamma^\mathsf{fg})^3}{(\gamma^\mathsf{fg})^3+1}.
\end{equation}
\end{corollary}

The greedy submodularity ratio can be obtained after the forward greedy algorithm is completed by analyzing $\mathcal{O}(\textstyle\binom{n}{K})$ inequalities. Since $\gamma^{\mathsf{fg}}\geq \gamma$, the performance guarantee in \eqref{eq: greedyversionPerformanceGuarantee} is better than~\eqref{eq: performanceguaranteeMatroidOpt}.
 Notice that $\gamma^{\mathsf{fg}}$ changes with the constraint set of the problem since the inequalities defining $\gamma^{\mathsf{fg}}$ would then be different. In contrast, submodularity ratio $\gamma$ depends only on the objective function.

Next, we substitute $f = -F_\epsilon$ and $\mathcal{F} = \tilde{\mathcal{C}}_K$ into the performance guarantee \eqref{eq: greedyversionPerformanceGuarantee} of the general setting~\eqref{eq:matroidconstrained optimization}. 
\looseness=-1
\begin{corollary}\label{cor:fwcor}
	  Suppose we apply Algorithm \ref{alg:ALG4} to Problem~\eqref{eq:mainproblem}. Denote the actuator set returned as $S^\mathsf{f}_\epsilon$ and the greedy submodularity ratio of $-F_\epsilon$ as $\gamma^{\mathsf{fg}}_\epsilon$. Then, $S^\mathsf{f}_\epsilon$ satisfies
\begin{equation}
  \frac{F_\epsilon(\emptyset)-F_\epsilon(S^\mathsf{f}_\epsilon)}{F_\epsilon(\emptyset)-F_\epsilon(S^*)}\geq\frac{(\gamma^\mathsf{fg}_\epsilon)^3}{(\gamma^\mathsf{fg}_\epsilon)^3+1}.
    \label{eq: boundwithEMPTYSET}
  \end{equation}
\end{corollary}

Since the forward greedy algorithm starts expanding from the empty set, performance guarantees can only assess $f(S^\mathsf{f})$ by considering $f(\emptyset)$ as the reference. If $f(\emptyset)=0$, the performance guarantee~\eqref{eq: greedyversionPerformanceGuarantee} is reduced to ${f(S^\mathsf{f})}/{f(S^*)}\geq {\gamma^3}/{(1+\gamma^3)}$. In this case, we only lose a fraction of the optimal objective by adopting the forward greedy algorithm. However, for our actuator placement problem $F_\epsilon(\emptyset) = n\epsilon^{-1}$, and the performance guarantee~\eqref{eq: boundwithEMPTYSET} is equivalent to
\begin{equation}
\label{eq:estimationOFFepsilon}
F_\epsilon(S_\epsilon^\mathsf{f})\leq \frac{1}{(\gamma^\mathsf{fg}_\epsilon)^3+1}F_\epsilon(\emptyset)+\frac{(\gamma^\mathsf{fg}_\epsilon)^3}{(\gamma^\mathsf{fg}_\epsilon)^3+1}F_\epsilon(S^*).
\end{equation}
Since $\epsilon$ is a small positive number and $n$ is in general large, the guarantee above can be loose.\footnote{If there exists an initial actuator set $S^\text{ini}\neq \emptyset$ rendering the system controllable, we can potentially mitigate this issue, since the reference of the guarantee would then be given by $F_\epsilon(S^\text{ini})$. Clearly, such applications also allow to set $\epsilon=0$, and drop structural controllability constraints.}
In the next section, we consider a variant of the greedy algorithm that comes along with a performance guarantee that does not depend on~$F_\epsilon(\emptyset)$. 

\section{Reverse Greedy Algorithm}
\label{sec:reverseGreedy}

To derive an alternative guarantee, we consider the reverse greedy algorithm (also called the stingy or greedy descent). This algorithm starts from the full set, and at each iteration, excludes the node with the least marginal gain from the actuator set of the previous iteration until a solution is reached. Such an approach allows to have the reference as $F_\epsilon(V)$, which is significantly smaller than $F_\epsilon(\emptyset)$ in practice. 

\subsection{Properties of the objective}
For the reverse greedy algorithm, we reformulate our metric as $F^\mathsf{r}_\epsilon(R):=F_\epsilon(V\setminus R),$ for all $R\subset V$.
The following lemma characterizes the properties of this function.

\begin{lemma}
\label{lem:weakrev}
The set function $F^\mathsf{r}_\epsilon$ is strictly increasing, weakly submodular with submodularity ratio $\gamma^\mathsf{r}_\epsilon>0$ and weakly supermodular with curvature $\alpha^\mathsf{r}_\epsilon<1$.
\end{lemma}
\begin{IEEEproof}
Regarding the strict monotonicity, suppose $S_1\subsetneqq S_2$. Since $F_\epsilon$ is strictly decreasing and $V\setminus S_2 \subsetneqq V\setminus S_1$, $F_\epsilon(V\setminus S_2)>F_\epsilon(V\setminus S_1)$, which implies $F^\mathsf{r}_\epsilon(S_2)>F^\mathsf{r}_\epsilon(S_1)$. Due to strict monotonicity of $F^\mathsf{r}_\epsilon$, it follows readily from the equalities shown in~\eqref{eq:extentionofGammaAlphaDef} that the submodularity ratio is strictly greater than 0 and the curvature is strictly less than 1. Thus, $F^\mathsf{r}_\epsilon$ is weakly submodular with $\gamma^\mathsf{r}_\epsilon>0$ and weakly supermodular with $\alpha^\mathsf{r}_\epsilon<1$. 
\end{IEEEproof}

Recall that the submodularity ratio of $-F_\epsilon$ is $\gamma_\epsilon^\mathsf{f}$. Now denote its curvature as $\alpha_\epsilon^\mathsf{f}$, which can easily be shown to satisfy  $\alpha_\epsilon^\mathsf{f}<1.$ The following connects $(\gamma_\epsilon^\mathsf{f},\alpha_\epsilon^\mathsf{f})$ and $(\gamma_\epsilon^\mathsf{r},\alpha_\epsilon^\mathsf{r})$.

\begin{proposition}\label{prop:rechar}
$\gamma_\epsilon^\mathsf{r} = 1-\alpha_\epsilon^\mathsf{f}$ and $\alpha_\epsilon^\mathsf{r} = 1-\gamma_\epsilon^\mathsf{f}.$
\end{proposition}

The proof is relegated to Appendix~\ref{app:props}. The proposition above provides an insight into how the submodularity ratio and the curvature of $F_\epsilon^\mathsf{r}$ relates to those of $-F_\epsilon$.

\subsection{Reformulation of the constraint set}
The reverse greedy algorithm has to return an exclusion set $R^\mathsf{r}$ such that the resulting actuator set $V\setminus R^\mathsf{r}$ contains $K$ nodes, and renders the system structurally controllable, that is, $V\setminus R^\mathsf{r}\in \mathcal{C}_K$. We collect all such exclusion sets and form ${\mathcal{R}}_{K}=\{R\,|\, V\setminus R\in {\mathcal{C}}_{K}\}$.  Suppose after the $t^\text{th}$ node exclusion of the reverse greedy algorithm, all the nodes excluded form a set $R^t=\{r_1,\ldots,r_t\}$, where $r_i\in V$ for all $i$. 
The set $R^t$ has to be a subset of some set in ${\mathcal{R}}_K$ for any $t = 1,\ldots,N$, where $N:=n-K$, since otherwise when $N$ exclusions are completed, the resulting actuator set would not belong to~$\mathcal{C}_K$. Thus, define
$\tilde{\mathcal{R}}_K:= \{Q\,|\, \exists R\in \mathcal{R}_K \text{ such that } Q\subset R \}$, and reformulate Problem \eqref{eq:mainproblem2} as
\begin{equation}
 \underset{R\subset V}{\min}\ F^\mathsf{r}_\epsilon(R) \ \ \mathrm{s.t.} \ R\in \tilde{\mathcal{R}}_{K} \text{ and $|R|= N$}.
\label{eq:mainproblemanotherform}
\end{equation}
The strict monotonicity of $F^\mathsf{r}_\epsilon(R)$ again ensures that the optimal solution to Problem \eqref{eq:mainproblemanotherform} coincides with that~of Problem~\eqref{eq:mainproblem2}. Note that cardinality constraint in \eqref{eq:mainproblemanotherform} can equivalently be replaced with an inequality constraint $|R|\geq N$.

 Next, we show that $\mathcal{\tilde{R}}_K$ characterizes a matroid.

\begin{proposition}
\label{prop: matroid2}
$\mathcal{M}^\mathsf{r}=(V,\mathcal{\tilde{R}}_K)$ is a matroid.
\end{proposition}

\begin{IEEEproof}
We prove that $(V,\mathcal{\tilde{R}}_K)$ is the dual of $(V,\mathcal{\tilde{C}}_K)$, that is, $\mathcal{\tilde{R}}_K=\mathcal{\tilde{C}}_K^*$. Note that we can then invoke Lemma~\ref{prop:dualmatroid} showing that the dual of a matroid is also a matroid. For any $Q\in \mathcal{\tilde{R}}_K$, according to the definition of $\mathcal{\tilde{R}}_K$,  there exists $R\in \mathcal{{R}}_K$ such that $Q\subset R$. We have $S = V\setminus R\in \mathcal{C}_K$ for $R\in \mathcal{{R}}_K$. Since $|S| = K$, $S$ is a maximum independent set in $\tilde{\mathcal{C}}_K$. Considering $Q\subset V\setminus S$, we conclude $Q\in\mathcal{\tilde{C}}^*_K.$ Conversely, for any $Q\in \mathcal{\tilde{C}}^*_K$, there exists a maximum independent set $S\in \mathcal{\tilde{C}}_K$ such that $Q\subset V\setminus S$. From the definition of $\mathcal{\tilde{C}}_K$, we know that $S\in \mathcal{{C}}_K$. Thus, we obtain $Q\in \mathcal{\tilde{R}}_K$. This concludes the equivalence of $\mathcal{\tilde{R}}_K$ and $\mathcal{\tilde{C}}_K^*$.
\end{IEEEproof}

Similar to the discussions in Section \ref{sec:Forwardgreedyalgorithm},
by restricting the iterates of the reverse greedy algorithm to lie in the set collection $\mathcal{\tilde{R}}_K$, we obtain a final exclusion set in ${\mathcal{R}}_K$ with cardinality $N$. This implies that the final actuator set lies in~$\mathcal{{C}}_K$.
\looseness=-1

\subsection{Performance guarantee}
In the previous section, we showed that the objective function $F^\mathsf{r}_\epsilon$ is  $\gamma_\epsilon^\mathsf{r}$-submodular and $\alpha_\epsilon^\mathsf{r}$-supermodular in Lemma~\ref{lem:weakrev}, and the feasible region of $\mathcal{\tilde{R}}_K$ characterizes a matroid in Proposition~\ref{prop: matroid2}. Thus, Problem~\eqref{eq:mainproblemanotherform} falls into the following class of optimization problems:
\begin{equation}
\begin{aligned}
& \underset{R\subset V}{\min}
& & f(R),\,\text{strictly increasing,}\text{ $\gamma$-submodular},\\ &&&\text{and $\alpha$-supermodular} \\
&\ \mathrm{s.t.} & & R\in \mathcal{F}\text{,\, $\mathcal{M}=(V,\mathcal{F})$ is a matroid,} |R|\geq N, 
\end{aligned}
\label{eq:generalized reverse version}
\end{equation}
where the cardinality of maximum independent sets in $\mathcal{F}$ is $N$.\footnote{The performance guarantee we derive in this section will be valid as long as the cardinality of maximum independent sets in $\mathcal{F}$ are larger than or equal to $N$, since this would ensure the feasibility of the problem.}
Let $R^*$ denote its optimal solution. Clearly, $R^*$ is the set complement of $S^*$, that is, $R^*=V\setminus S^*.$

Define set function $f^\mathsf{o}$ such that $f^\mathsf{o}(R)=f(V\setminus R)$ for all $R$. In Problem~\eqref{eq:mainproblemanotherform}, $f^\mathsf{o}$ corresponds to $F_\epsilon$.
Observe that the forward greedy algorithm applied to the minimization of the function~$f$ is equivalent to the reverse greedy algorithm applied to the minimization of the function $f^\mathsf{o}.$ This algorithm is presented in Algorithm~\ref{alg:ALGreverseGreedy}. Different from Algorithm \ref{alg:ALG4}, at each iteration, Algorithm~\ref{alg:ALGreverseGreedy} implements the feasibility check on the node with the least marginal gain. 

For Algorithm~\ref{alg:ALGreverseGreedy}, the following definitions are in order. We define $\rho_j(R) := f(R\cup j)-f(R)$, $\rho_{t} := f(R^{t})-f(R^{t-1})$ and $r_t := R^t\setminus R^{t-1}$. The set $U^t$ denotes the set of nodes having been considered by the feasibility check before $r_{t+1}$. The final exclusion set is $R^\mathsf{r}:=R^N$, and it lies in ${\mathcal{R}}_K$.
\begin{algorithm}[t]
        \caption{Reverse Greedy Algorithm over Matroid}
        \label{alg:ALGreverseGreedy}
        \begin{algorithmic}
        \Require set function $f^\mathsf{o}$, ground set $V$, matroid $(V,\mathcal{F})$
        \Ensure exclusion set $R^\mathsf{r}$
        \Function {ReverseOverMatroid}{$f^\mathsf{o}, V,{\mathcal{F}}$}
        
        \State {$R^0 = \emptyset$, $U^0 = \emptyset$, $t=1$}
        \While{$U^{t-1}\neq V$ { and } $|R^{t-1}|<N$}
            \State ${j^*(t)}= \arg\min_{j\in V\setminus U^{t-1}}\rho_j(R^{t-1})$
            \If {$R^{t-1}\cup j^*(t)\notin\mathcal{\mathcal{F}}$}
        \State $U^{t-1}\gets U^{t-1}\cup j^*(t)$
        \Else
        \State $\rho_{t} \gets\rho_{j^*(t)}(R^{t-1})$ and $r_t = j^*(t)$
        \State $R^{t} \gets R^{t-1}\cup j^*(t)$ and $U^{t}\gets U^{t-1}\cup j^*(t)$
        \State $t\gets t+1$
    \EndIf    
            \EndWhile
            \State $R^\mathsf{r} \gets R^{t-1} $
        \EndFunction
        \end{algorithmic}
    \end{algorithm}

A special case of Problem~\eqref{eq:generalized reverse version} was previously shown to be hard to approximate. Specifically, for the problem of minimizing a submodular increasing function over only a cardinality lower bound, the work in~\cite{svitkina2011submodular} shows that there is no bicriteria approximation performing better than $o(\sqrt{n/\text{log}\, n})$, where $n$ is the cardinality of the ground set.\footnote{Bicriteria approximation refers to approximating both the constraint requirement and the optimal objective. We refer to \cite{svitkina2011submodular} for the exact description.} 
Next, we extend this result by providing novel counterexamples showing that a strictly positive submodularity ratio and a curvature bounded away from $1$ is indispensable to obtain any meaningful performance guarantee for Problem~\eqref{eq:generalized reverse version}. The proofs of Propositions~\ref{prop: impossibilityResults1} and~\ref{prop: impossibilityResults2} are relegated to Appendix~\ref{app:props}.

\begin{proposition}
\label{prop: impossibilityResults1}
In Problem \eqref{eq:generalized reverse version}, one cannot derive any upper bound on $(f(R^\mathsf{r})-f(\emptyset))/(f(R^*)-f(\emptyset))$ if no strictly positive lower bound on $\gamma$ is known. 
\end{proposition}
\begin{proposition}
\label{prop: impossibilityResults2}
In Problem \eqref{eq:generalized reverse version}, one cannot derive any upper bound less than $N$ on $(f(R^\mathsf{r})-f(\emptyset))/(f(R^*)-f(\emptyset))$, if no upper bound less than $1$ is known for $\alpha$.
\end{proposition}

The propositions above conclude that we have to utilize both the submodularity ratio and the curvature. 
Our second main result is shown in the following theorem.

\begin{theorem}
\label{thm:reverseGuarantee}
If Algorithm \ref{alg:ALGreverseGreedy} is applied to \eqref{eq:generalized reverse version}, then
\begin{equation}
\label{eq: reverseBound}
\frac{f(R^\mathsf{r})-f(\emptyset)}{f(R^*)-f(\emptyset)} \leq  \frac{\gamma}{1-\gamma}\big((2N+1)^{\frac{1-\gamma}{\gamma(1-\alpha)}}-1\big).
\end{equation}
\end{theorem}

The proof extends the linear programming proofs utilized by \cite{conforti1984submodular}, which considers the maximization of increasing submodular functions over matroid constraints, and by \cite{bian2017guarantees}, which considers the maximization of increasing, nonsubmodular nonsupermodular functions over cardinality constraints. In contrast, our proof applies to the minimization of increasing, nonsubmodular nonsupermodular functions over matroids.

The main idea of the proof is to provide a series of inequalities that upperbound $\rho_t$ by $f(R^*)-f(\emptyset)$, for each iteration $t$. This way, $f(R^\mathsf{r})-f(\emptyset)=\sum_{t=1}^N\rho_t$ has an upper bound expressed by $f(R^*)$. For the following lemma, we recall that $R^t:=\{r_1,\ldots,r_t\}$ is the set obtained by the greedy algorithm after the exclusion~of~$r_t$.

\begin{lemma}
For any $t\in\{0,\ldots,N-1\}$, $\rho_t$ satisfies
\begin{equation}
\begin{aligned}
f(R^*)-f(\emptyset)\geq & (1-\frac{1}{\gamma})\sum_{i:r_i\in R^t\setminus R^*}\rho_i + \sum_{i:r_i\in R^t\cap R^*}\rho_i 
\\
&+ (1-\alpha)(N-t)\rho_{t+1}.
\end{aligned}
\label{eq:LMM1}
\end{equation}
\end{lemma}

\begin{IEEEproof}
Suppose $R^* = \{r^*_1,\ldots, r^*_{N}\}$. Rewrite $f(R^*\cup R^t)$ as two telescoping sums
$f(R^*\cup R^t) = f(R^*)+\sum_{i=1}^t\rho_{r_i}(R^* \cup R^{i-1}),$ and $f(R^*\cup R^t)=f(R^t)+\sum_{k=1}^{N}\rho_{r^*_k}(\{r^*_1,\ldots ,r^*_{k-1}\}\cup R^t),$ {which is directly obtained from the definition of $\rho_{r}$.} 
For any $i$ such that $r_i\in R^*\cap R^t $, we have $\rho_{r_i}(R^* \cup R^{i-1}) = 0$. Using this, and the fact that both telescoping sums above are equal to $f(R^*\cup R^t)$, we obtain
\begin{equation}
\begin{aligned}
f(R^*)+\sum_{i:r_i\in R^t\setminus R^*}\rho_{r_i}(R^* \cup R^{i-1})&=f(R^*\cup R^t)\\
&\hspace{-4cm}=f(R^t)+\sum_{k=1}^{N}\rho_{r^*_k}(\{r^*_1,\ldots,r^*_{k-1}\}\cup R^t).
\end{aligned}
\label{eq:baseequality}
\end{equation}
Invoking the definitions of submodularity ratio and curvature, for each $i$ such that $r_i\in R^t\setminus R^*$,  we have
\begin{equation}
\begin{aligned}
\rho_{r_i}(R^*\cup R^{i-1})\leq \frac{1}{\gamma}\rho_i,
\end{aligned}
\label{eq:leftinequality}
\end{equation}
and for any $k \in \{1,\ldots, N\}$,
\begin{equation}
\begin{aligned}
\rho_{r^*_k}(\{r^*_1,\ldots, r^*_{k-1}\}\cup R^t)\geq (1-\alpha )\rho_{r^*_k}(R^t).
\end{aligned}
\label{eq:rightinequality}
\end{equation}
    By the definition of a matroid, there exists $R_c^t=\{r^*_{c_1},\ldots,r^*_{c_{N-t}}\}\subset R^*\setminus R^t$ such that $R_c^t\cup R^t\in \mathcal{{F}}$. Consequently, for any $1\leq i\leq N-t$, we have $R^t\cup r^*_{c_i}\in \mathcal{{F}}$. Thus, adding $r^*_{c_i}$ to $R^{t}$ has to be feasible in the matroid. If $\rho_{r^*_{c_i}}(R^{t})<\rho_{r_{t+1}}(R^{t})=\rho_{t+1}$, $r^*_{c_i}$ could be added to $R^{t}$ to form $R^{t+1}$ instead of $r_{t+1}$. This yields a contradiction, implying that the inequality $\rho_{r^*_{c_i}}(R^{t})\geq\rho_{r_{t+1}}(R^{t})=\rho_{t+1}$ holds for any $i$. Hence, we obtain
\begin{equation}
\begin{aligned}
\sum_{k=1}^N\rho_{r^*_k}( R^t)\geq \sum_{r^*_k\in R^t_c}\rho_{r^*_k}(R^t)\geq (N-t)\rho_{t+1}.
\end{aligned}
\label{eq:rightinequality2}
\end{equation}
Next, by substituting (\ref{eq:leftinequality}), (\ref{eq:rightinequality}) and (\ref{eq:rightinequality2}) into (\ref{eq:baseequality}), we obtain
{\medmuskip=.4mu\thinmuskip=.4mu\thickmuskip=.4mu\begin{equation}
\begin{aligned}
f(R^*)-f(\emptyset)+\sum_{\small\substack{i:\\r_i\in R^t\setminus R^*}}\frac{\rho_{i}}{\gamma}\geq \sum_{i=1}^t\rho_{i} +(1-\alpha)(N-t)\rho_{t+1}.
\end{aligned}
\label{eq:BEFOREfinal}
\end{equation}}By grouping the terms in (\ref{eq:BEFOREfinal}), we obtain (\ref{eq:LMM1}).
\end{IEEEproof}

Next, we can construct a linear program where the  solution provides an upper bound for $\frac{f(R^\mathsf{r})-f(\emptyset)}{f(R^*)-f(\emptyset)}$.

\begin{IEEEproof}[Proof of Theorem~\ref{thm:reverseGuarantee}] Let $x_i = \rho_{i}/(f(R^*)-f(\emptyset))$, we have $(f(R^\mathsf{r})-f(\emptyset))/(f(R^*)-f(\emptyset)) = \sum_{i=1}^{N}x_i.$ Note that $x_i\geq 0$ for all $i$. Suppose $R^*\cap R^\mathsf{r} = \{r_{i_1},r_{i_2},\ldots \}$. To give an upper bound for this ratio, we exploit the inequalities \eqref{eq:LMM1} and build the following linear programming problem to compute the largest possible sum, $\sum_{i=1}^{N}x_i$,
{\medmuskip=.5mu\thinmuskip=.5mu\thickmuskip=.5mu\begin{equation}
\begin{aligned}
&Z(N,\gamma,\alpha)=\max \sum_{i=1}^{N}x_i,\text{ }{\rm { s.t. }}\text{ } x_i\geq 0 \text{ }{\rm { and }}\text{ }\\
&\scalemath{0.94}{{\scriptsize	\begin{bmatrix}
    (1-\alpha)N &  &  &  &    &\\
    -(1-\gamma)/\gamma & \substack{(1-\alpha)\\\times(N-1)} &  &   &    &\\
    \vdots &  \vdots & \ddots &  &    & \\
    -(1-\gamma)/\gamma & -(1-\gamma)/\gamma & \dots & \substack{(1-\alpha)\\\times(N-i_1+1)}   &    & \\
    -(1-\gamma)/\gamma & -(1-\gamma)/\gamma& \dots & 1 & \ddots& & \\
    \vdots & \vdots & \dots & \vdots  & \vdots  & \\
    -(1-\gamma)/\gamma & -(1-\gamma)/\gamma & \dots & 1  & \dots &  (1-\alpha)
\end{bmatrix}\hspace{-.1cm}
\begin{bmatrix}
    x_1\\
    x_2\\
    \vdots \\
    x_{i_1} \\
    \vdots \\
    x_N
\end{bmatrix}\hspace{-.1cm}
\leq
\begin{bmatrix}
    1\\
    1\\
    \vdots \\
    1 \\
    1 \\
    \vdots \\
    1
\end{bmatrix}.}}
\end{aligned}
\label{eq:rlp}
\end{equation}}To get an upper bound for $\sum_{i=1}^N x_i$, we consider the following relaxed problem where the unit entries in \eqref{eq:rlp} are replaced by $-(1-\gamma)/\gamma$,
{\medmuskip=.58mu\thinmuskip=.58mu\thickmuskip=.58mu\begin{equation}
\begin{aligned}
&\bar{Z}(N,\gamma,\alpha)=\max \sum_{i=1}^{N}x_i,\text{ }{\rm { s.t. }}\text{ } x_i\geq 0 \text{ }{\rm { and }}\text{ }\\
&{\scriptsize\begin{bmatrix}
    (1-\alpha)N &  &  & \\
    -(1-\gamma)/\gamma & (1-\alpha)(N-1) &  &   \\
    \vdots &  \vdots & \ddots & \\
    -(1-\gamma)/\gamma & -(1-\gamma)/\gamma & \dots & (1-\alpha)    \\

\end{bmatrix}
\begin{bmatrix}
    x_1\\
    x_2\\
    \vdots \\
    x_N
\end{bmatrix}
\leq
\begin{bmatrix}
    1\\
    1\\
    \vdots \\
    1 \\
\end{bmatrix}.}
\end{aligned}
\label{eq:lp2}
\end{equation}}Since we require that $x_i\geq0$ for any $i$, any feasible solution to~\eqref{eq:lp2} is also feasible to~\eqref{eq:rlp}. Thus, $\bar{Z}(N,\gamma,\alpha)\geq Z(N,\gamma,\alpha)$ for any $N$, $\gamma$ and $\alpha$. We claim that the optimum $x^*$ of Problem~\eqref{eq:lp2} makes all the inequality constraints tight. This is easily seen by rewriting the inequalities as $x_t\leq\frac{1}{(1-\alpha)(N-t+1)}(1+\frac{1-\gamma}{\gamma}\sum^{t-1}_{i=1}x_i),\ t=1,\ldots,N.$
Notice that for submodular functions, we have $\gamma=1$. Using the claim above and considering the fact that $x^{-1}<\text{ln}(x+1/2)-\text{ln}(x-1/2)$ for any $x\geq 1$, we can directly obtain the following guarantee
{\medmuskip=.66mu\thinmuskip=.66mu\thickmuskip=.66mu\begin{equation}
\dfrac{(f(R^\mathsf{r})-f(\emptyset))}{(f(R^*)-f(\emptyset))}\leq\sum_{i=1}^{N}x_i^*=\sum_{i=1}^{N}\dfrac{1}{i(1-\alpha)}\leq\dfrac{\text{ln} (2N+1)}{1-\alpha}.
\label{eq:reverseguaranteeForSubmodularFunction}
\end{equation}}

Next, we focus our efforts on the case in which $0<\gamma<1$ and $0\leq\alpha<1$. We obtain
\begin{equation}\label{eq:actual_rev_guar}\bar{Z}(N,\gamma,\alpha) = \frac{1}{b}\prod_{i=1}^{N}\bigg(1+\frac{b}{(N-i+1)(1-\alpha)}\bigg)-\frac{1}{b},\end{equation}
where $b=(1-\gamma)/\gamma.$
Considering $\bar{Z}(N,\gamma,\alpha)\geq Z(N,\gamma,\alpha)\geq (f(R^\mathsf{r})-f(\emptyset))/(f(R^*)-f(\emptyset))$ and 
{\medmuskip=.6mu\thinmuskip=.6mu\thickmuskip=.6mu\begin{equation}
\begin{aligned}
\bar{Z}(N,\gamma,\alpha)  &= b^{-1}\text{exp}\bigg(\sum_{i=1}^{N} \text{ln}\Big( 1+\frac{b}{(N-i+1)(1-\alpha)}\Big)\bigg)-b^{-1}\\
& \leq b^{-1}\text{exp}\bigg(\frac{b}{(1-\alpha)}\sum_{i=1}^{N}\frac{1}{(N-i+1)}\bigg)-b^{-1}\\
& \leq b^{-1}\text{exp}\bigg(\frac{b}{(1-\alpha)}\Big(\text{ln}(N+\frac{1}{2})-\text{ln}\frac{1}{2}\Big)\bigg)-b^{-1}
\\
& = b^{-1}(2N+1)^{\frac{b}{(1-\alpha)}}-b^{-1}. 
\end{aligned}
\label{eq:RbarEstimation}
\end{equation}}The first equality rewrites the multiplication in $\bar{Z}$ into an exponential sum. The inequalities follow from the fact that $\text{ln}(1+x)<x$ for any $x>0$ and $x^{-1}<\text{ln}(x+1/2)-\text{ln}(x-1/2)$ for any $x\geq 1$. By substituting $b=(1-\gamma)/\gamma$ back into the last term in \eqref{eq:RbarEstimation}, we get \eqref{eq: reverseBound}.
\end{IEEEproof}

Let $Z_u(N,\gamma,\alpha)= \frac{\gamma}{1-\gamma}\big((2N+1)^{\frac{1-\gamma}{\gamma(1-\alpha)}}-1\big)$. Table~\ref{tab: ComparisonW} illustrates how well the upper bound $Z_u$ approximates the original guarantee $\bar{Z}$, stated in \eqref{eq:actual_rev_guar}. For the supermodular case $\alpha=0$, we obtain the guarantee $ \frac{\gamma}{1-\gamma}\big((2N+1)^{\frac{1-\gamma}{\gamma}}-1\big).$ Via the upperbound in~\eqref{eq:reverseguaranteeForSubmodularFunction}, for the submodular case $\gamma=1$, we obtain $\frac{\text{ln} (2N+1)}{1-\alpha}.$
As a remark, we can verify that the guarantee in~\eqref{eq:actual_rev_guar} is not tight.
Suppose $f$ is modular, that is, both supermodular and submodular. Then, we have
$\frac{(f(R^\mathsf{r})-f(\emptyset))}{(f(R^*)-f(\emptyset))}\leq\text{ln} (2N+1).$
However, modularity of $f$ implies that the greedy algorithm returns the optimal solution~\cite{edmonds1971matroids}. One reason for this looseness is that, to ensure the tightness of the relaxation from~\eqref{eq:rlp} to~\eqref{eq:lp2}, we must have $R^*\cap R^t=\emptyset$, which then contradicts the modularity of the objective function. To the best of our knowledge, Theorem~\ref{thm:reverseGuarantee} provides the first performance guarantee for the reverse greedy algorithm for this setting involving the submodularity ratio and the curvature.\footnote{\cite[Thm 7]{sviridenko2017optimal} offers a guarantee for the forward greedy algorithm applied to minimizing increasing functions over a matroid as in Problem~\eqref{eq:generalized reverse version}. This is given by $1/(1-c)$, where $c$ quantifies how far a function is from being modular. This novel notion is a significantly stronger requirement than having both the submodularity ratio and the curvature simultaneously, see~\cite[(6)]{sviridenko2017optimal}. Hence, it is not possible to compare it with our guarantee other than the case of a modular objective. In that case, setting $c=0$ confirms the optimality of the greedy algorithm. Note that computing this novel notion requires an exhaustive enumeration and it does not allow any greedy computation, which can limit its applications.}

\renewcommand{\arraystretch}{1.3}
\begin{table}[t]
\centering
\caption{Comparison between performance guarantees $\bar{Z}$ and $Z_u$}
\resizebox{.2\textwidth}{!}{\begin{tabular}{|c|c|c|}
\hline
 $(N,\gamma,\alpha)$&     $\bar{Z}$ &   $Z_u$   \\ \hline
 (20,0.9,0.1)&  4.87 & 5.25  \\ \hline
 (100,0.9,0.1)& 7.87 & 8.32 \\ \hline
 (20,0.99,0.1)& 4.07 & 4.21 \\ \cline{1-3}
\end{tabular}}
\label{tab: ComparisonW}
\end{table}
\renewcommand{\arraystretch}{1}

Similar to our analysis in Section \ref{sec:Forwardgreedyalgorithm}, 
we propose computationally more efficient approaches to deriving both the submodularity ratio and the curvature.

\begin{corollary}Let $\gamma^\mathsf{rg}$ be the largest $\hat{\gamma}$ that satisfies
$\rho_{r_t}(R\cup R^{t-1})\leq \hat{\gamma}^{-1} \rho_t,$ for all $t\leq N$, and $R$ with $|R|=N$. Let $\alpha^\mathsf{rg}$ be the smallest $\hat{\alpha}$ that satisfies $\rho_{r}(R\cup R^t)\geq (1-\hat{\alpha} )\rho_{r^*_k}(R^t),$
for all $t\leq N$, $R $ with $|R| = N-1$. Then, $\gamma^{\mathsf{rg}}$ is called the greedy submodularity ratio for the reverse greedy algorithm, with $\gamma^{\mathsf{rg}}\geq \gamma$, and $\alpha^{\mathsf{rg}}$ is called the greedy curvature for the reverse greedy algorithm, with $\alpha^{\mathsf{rg}}\leq \alpha$. The performance guarantee is given by
\begin{equation}
\label{eq: greedyversionPerformanceGuarantee2}
\frac{f(R^\mathsf{r})-f(\emptyset)}{f(R^*)-f(\emptyset)}\leq Z_u(N,\gamma^\mathsf{rg},\alpha^\mathsf{rg}).
\end{equation}
\end{corollary}

The greedy submodularity ratio above can be obtained after the reverse greedy algorithm is completed by analyzing $N\textstyle\binom{n}{N}$ inequalities, whereas the greedy curvature can be obtained by analyzing $N\textstyle\binom{n}{N-1}$ inequalities. Since $\gamma^{\mathsf{rg}}\geq \gamma$ and $\alpha^{\mathsf{rg}}\leq \alpha$, it can easily be verified that $Z_u(N,\gamma^\mathsf{rg},\alpha^\mathsf{rg})\leq Z_u(N,\gamma,\alpha).$

Substitute $f=F^\mathsf{r}_\epsilon$ and $\mathcal{F}=\tilde{\mathcal{R}}_K$ to conclude the following.

\begin{corollary}\label{cor:revcor}
Suppose we apply Algorithm~\ref{alg:ALGreverseGreedy} to Problem~\eqref{eq:mainproblemanotherform}. Denote the exclusion set returned as $R^\mathsf{r}_\epsilon$ and the greedy submodularity ratio of $F^\mathsf{r}_\epsilon$ as $\gamma_\epsilon^\mathsf{rg}$ and the greedy curvature of $F^\mathsf{r}_\epsilon$ as $\alpha_\epsilon^\mathsf{rg}$. Then, $R^\mathsf{r}_\epsilon$ satisfies
{\medmuskip=1mu\thinmuskip=1mu\thickmuskip=1mu\begin{equation}
\begin{aligned}
\label{eq: performanceguaranteeforreverse}
&\dfrac{F_\epsilon(V\setminus R_\epsilon^\mathsf{r})-F_\epsilon(V)}{F_\epsilon(S^*)-F_\epsilon(V)}\leq Z_u^\mathsf{rg},\ \text{or\ equivalently,}\\
   &F_\epsilon(V\setminus R_\epsilon^\mathsf{r})\leq Z_u^\mathsf{rg} {F_\epsilon(S^*)}+(1-Z_u^\mathsf{rg} )F_\epsilon(V),
\end{aligned}\end{equation}}where $Z_u^\mathsf{rg} := Z_u(N,\gamma^\mathsf{rg},\alpha^\mathsf{rg})$.
\end{corollary}

{In contrast to the forward greedy guarantee, $F_\epsilon(\emptyset)$ does not appear in the guarantee above, which is generally large. On the other hand, the guarantee above scales with the problem size, specifically, with $N=n-K.$} {In the numerics, we show that both greedy algorithms achieve comparable performance in our problem, and at the same time much better performance than what the theoretical guarantees suggest. In practice, it could be useful to implement both greedy algorithms (which can be done efficiently with polynomial time complexity) and choose the best out of the two.}

\section{Implementation Aspects}
\label{sec:FeasibilityCheck}

Two issues have to be addressed to implement Algorithms~\ref{alg:ALG4} and~\ref{alg:ALGreverseGreedy} for the actuator placement problem. First, we have to select a metric-modifying parameter $\epsilon$. Second, we need feasibility check methods for the set collections $\mathcal{\tilde{C}}_K$ and $\mathcal{\tilde{R}}_K$.\
\subsection{An algorithm for picking a metric-modifying parameter $\epsilon$}
\label{sec: algorithmforproperep}

The performance guarantees \eqref{eq:estimationOFFepsilon} and \eqref{eq: performanceguaranteeforreverse} relate to $F_\epsilon$ instead of the original metric $F$. {On the one hand, if $\epsilon$ is large, a performance guarantee on $F_\epsilon$ may not be applicable as a performance guarantee on $F$ since $F_\epsilon(S_\epsilon)<F(S_\epsilon)$.} On the other hand, if $\epsilon$ is small, the matrix $W_T(S)+\epsilon I$ may be close to singularity. Such ill-conditioned matrices can occur especially at the early stages of the forward greedy algorithm.\footnote{If there exists an initial actuator set $S^\text{ini}\neq \emptyset$ rendering the system controllable, invertibility is guaranteed without $\epsilon$. The reverse greedy also mitigates this issue, since it evaluates $F_\epsilon$ at structurally controllable systems. }
Denote the actuator set returned by a greedy algorithm applied to $F_\epsilon$ as $S_\epsilon$. This could be the solution returned by either Algorithm~\ref{alg:ALG4} or Algorithm \ref{alg:ALGreverseGreedy}.  Given an approximation factor $\xi>0$ as a design parameter, we propose an algorithm to pick $\epsilon$ such that $F(S_\epsilon)<(1+\xi)F_\epsilon(S_\epsilon)$. This inequality implies that guarantees in \eqref{eq:estimationOFFepsilon} and \eqref{eq: performanceguaranteeforreverse} translate into guarantees for the metric $F$.
The method is presented in Algorithm \ref{alg:epsilon}. Denote the eigenvalues of the controllability Gramian $W_T(S)$ as $\lambda_1(S)\leq\ldots \leq\lambda_n(S)$.\begin{algorithm}[t]
        \caption{Finding $\epsilon$ with provable performance}
        \label{alg:epsilon}
        \begin{algorithmic}
        \Require approximation factor $\xi$, initial value $\epsilon_0$
        \Ensure parameter $\epsilon$
        \Function {ProperEP}{$\xi,\epsilon_0$}
        \State {$i = 0$}
        \While{$\epsilon_i \geq \xi\lambda_1(S_{\epsilon_i})$}
        \State let $\epsilon_{i+1} \gets \frac{1}{2}\xi \lambda_1(S_{\epsilon_{i}})$ and $i\gets i+1$
            \EndWhile
        \EndFunction
        \end{algorithmic}
    \end{algorithm}\begin{proposition}\label{prop:epsilonalgproof}Suppose given any $\epsilon>0$, $W_T(S_\epsilon)$ is invertible. Then, for any approximation factor $\xi>0$ and any initial value $\epsilon_0>0$, Algorithm~\ref{alg:epsilon} returns $(\epsilon,S_\epsilon)$ pair that satisfies $F(S_\epsilon)<(1+\xi)F_\epsilon(S_\epsilon).$
\end{proposition}
\begin{IEEEproof}
If the controllability Gramian $W_T(S_\epsilon)$ is invertible, then we have that $\lambda_1(S_\epsilon)>0$. Since there are finitely many combinations of actuators, the set $\{\lambda_1(S_\epsilon)|\forall\epsilon>0\}$ has a positive lower bound, denoted as $\lambda_0$.
In the iterations of Algorithm~\ref{alg:epsilon}, it holds that $\epsilon_{i+1}<\frac{1}{2}\epsilon_{i}$, because $\epsilon_i >\xi\lambda_1(S_{\epsilon_i})$ and $\epsilon_{i+1} =\frac{1}{2}\xi\lambda_1(S_{\epsilon_i})$. Hence, there exists some $j$ such that $\epsilon_{j}<\xi \lambda_0\leq \xi\lambda_1(S_{\epsilon_{j}})$. Then, we obtain $F(S_{\epsilon_{j}})= \sum^n_{i=1}\frac{1}{\lambda_i(S_{\epsilon_{j}})} <\sum^n_{i=1}\frac{1+\xi}{\lambda_i(S_{\epsilon_{j}})+\epsilon_{j}}= (1+\xi) F_{\epsilon_{j}}(S_{\epsilon_{j}}).$ This inequality concludes the proof.\end{IEEEproof}

For the proof above, we assumed that given any $\epsilon>0$, the controllability Gramian $W_T(S_\epsilon)$ is invertible. This is a strong assumption since, as previously mentioned, structural controllability does not imply controllability.
In the numerics, we always ended up with a controllable system with any of the greedy algorithms. This can be explained either by the objective of the problem which is to minimize the average energy consumption or the choice of a large cardinality $K$. 

We now provide the resulting performance guarantees.
\begin{corollary}
Given the factor $\xi$, suppose we apply Algorithm~\ref{alg:epsilon} to pick $\epsilon$. From Corollaries~\ref{cor:fwcor} and~\ref{cor:revcor}, we have 
\begin{equation*}
\begin{aligned}
&F(S_\epsilon^\mathsf{f})<(1+\xi)\Big[ \frac{1}{(\gamma^\mathsf{fg}_\epsilon)^3+1}F(\emptyset)+\frac{(\gamma^\mathsf{fg}_\epsilon)^3}{(\gamma^\mathsf{fg}_\epsilon)^3+1}F(S^*)\Big],\\
&F(V\setminus R_\epsilon^\mathsf{r})<(1+\xi)\left[ Z_u^\mathsf{rg} {F(S^*)}+(1-Z_u^\mathsf{rg} )F_\epsilon(V)\right],\\
\end{aligned}
\end{equation*}
where $S^*$ is the optimal solution to \eqref{eq:mainproblem2}.
\end{corollary}

\subsection{Feasibility check over $\mathcal{\tilde{C}}_K$}
\label{sec:feasibilitycheck}
When applied to Problem \eqref{eq:mainproblem}, the forward greedy algorithm has to ensure that the actuator set returned by each iteration lies in $\tilde{\mathcal{C}}_K$. {The work of~\cite{liu2011controllability,Liu2011ControllabilityOC} proposes a method to determine whether a given set~$S$ with $|S| = K$ belongs to $\mathcal{C}_K$. This result is not directly applicable to answer whether an actuator set $S$ with $|S|<K$ returned by a greedy iteration belongs to $\mathcal{\tilde{C}}_K$. 
In the following, we extend the work of~\cite{liu2011controllability,Liu2011ControllabilityOC} for a feasibility check over $\mathcal{\tilde{C}}_K$ by constructing auxiliary bipartite graphs associating this check with the cardinality of a maximum matching and by formulating a maximum flow problem.}
\looseness=-1

We introduce the concept of matchings and bipartite graphs. An undirected graph is called bipartite and denoted as $(V^1, V^2, \mathsf{E})$ if its vertices are partitioned into $V^1$ and $V^2$ while any edge in $\mathsf{E}$ connects a vertex in $V^1$ to another in~$V^2$. A matching $m$ is a subset of $\mathsf{E}$ if no two edges in $m$ share a vertex in common. Given a subset $L$ of $V^1\cup V^2$, we say $L$ is covered by $m$ if any $v\in L$ is connected to an edge in $m$. Matching $m$ is maximum if it has the largest cardinality among all the matchings and is perfect if $V^2$ is covered.

Given the graph $G=(V,E)$ describing system \eqref{eq: systemmodel}, we first build the following auxiliary bipartite graph to determine whether an actuator set renders the system structurally controllable. Node sets $V'=\{v'_1,\ldots, v'_n\}$ and ${V''}=\{v''_1,\ldots, v''_n\}$ are built as two copies of $V =\{{v}_1,\ldots, {v}_n\}$. For any set $S\subset V$, two subsets $S'\subset V'$ and $S''\subset V''$ denote two copies of the set $S$. For the bipartite graph, we then have $V^1=V\cup S''$ and $V^2=V'$. Next, we define the edge sets. The set $\mathsf{E}$ consists of undirected edges connecting $v_{i}$ with $v'_{j}$ if $(v_i,v_j)\in E$, whereas the edge set $\mathsf{E}_1$ consists of undirected edges connecting $v'_k$ with $v''_k$ if $v_k\in S$. The bipartite graph is then defined by $\mathcal{H}_b(S) = (V\cup S'', V', \mathsf{E}\cup \mathsf{E}_1)$. 

{If the graph $G$ and the set~$S$ pair satisfies the accessibility condition}\footnote{\label{foot:9}{The graph $G$ and the set~$S$ pair satisfies the accessibility condition if for every node in~$G$ there is at least one directed path reaching that node from some node in~$S$~\cite{Liu2011ControllabilityOC}~\cite[Def.~2]{clark2017submodularity}. When $G$ is strongly connected, this condition is attained irrespective of the nodes chosen in $S$~\cite{clark2017submodularity}. We can invoke \cite[Thm~2]{Liu2011ControllabilityOC} because of the equivalence of structural controllability in leader selection and actuator placement we established in the proof of Proposition~\ref{prop: matroid}.}} {(which is implied by the strong connectivity assumption on~$G$), the set $S$ achieves structural controllability if and only if there exists a perfect matching in $\mathcal{H}_b(S)$, see~\cite[Thm~2]{Liu2011ControllabilityOC}.}
This equivalence directly follows from Hall's marriage theorem, which shows that there exists a perfect matching in $\mathcal{H}_b(S)$ if and only if, for any $U\subset V'$, the nodes in $U$ have at least $|U|$ unique in-neighbors~\cite{clark2017submodularity}.
Intuitively, to control any node, we would influence the states of its in-neighbors in the graph. Then, to steer the nodes in $U$ arbitrarily, this theorem implies that we should have at least $|U|$ in-neighbors. Otherwise, suppose two nodes share only a single in-neighbor. Then, these nodes would always be receiving a proportional influence, making it impossible to steer the system states arbitrarily. Using this result, \cite{clark2012leader} develops a feasibility check for leader selection. This method states that $S$ lies in $\tilde{\mathcal{C}}_K$ if and only if there is a maximum matching for the bipartite graph $\mathcal{H}_b(\emptyset)$ with all the nodes in the set $S'\subset V'$ unmatched. However, this statement is true only if we consider the minimum required cardinality for the structural controllability of the system, see the proof of \cite[Lemma~3]{clark2012leader}. Later in this section, we provide a counterexample where the feasibility check of \cite{clark2012leader} does not work.

We now provide our feasibility check in the following.

\begin{proposition}
\label{prop: feasibility}
Given the graph $G$, the cardinality limit $K$ and an actuator set $S$ with $|S|=k\leq K$, we have $S\in \mathcal{\tilde{C}}_K$ if and only if $|\bar{m}(S)| \geq n-K+k$, where $\bar{m}(S)$ is a maximum matching in $\mathcal{H}_b(S)$.
\end{proposition}

\begin{IEEEproof}
``$\Rightarrow$": If $S\in \mathcal{\tilde{C}}_K$, there exists $Q \in \mathcal{C}_K$ such that $S\subset Q$. We now invoke the equivalence result from~\cite[Thm~2]{Liu2011ControllabilityOC}. This implies the following. By finding a maximum matching $m$ in $\mathcal{H}_b(Q)$ that completely covers ${Q}''$ and then excluding from $m$ the edges incident with ${Q''}\setminus {S''}$, we can obtain a matching in $\mathcal{H}_b( S)$ containing $n-K+k$ edges. 

``$\Leftarrow$": We pick any maximum matching in $\mathcal{H}_b(S)$ and denote it as $m^*$. Suppose $P'$ is the largest subset in $V'$ whose elements are all missed by $m^*$, we know $|P'|\leq K-k$. Denote the edge subset $\mathsf{E}_P\subset \mathsf{E}$ as the set that contains all undirected edges adjacent to $v'_k$ for any $k$ such that $v_k \in P$. Clearly, $\mathsf{E}_P$ covers $P'$ and $m^*\cup \mathsf{E}_P$ covers $V'$. Since matching $m^*$ and matching $\mathsf{E}_P$ have no common vertices, $m^*\cup \mathsf{E}_P$ is a perfect matching in $\mathcal{H}_b(S\cup P)$, which means with the actuator set $ S\cup P$ the system is structurally controllable. Also considering $|S \cup P|\leq K$, and $S\cup P\in \mathcal{\tilde{C}}_K$, we obtain $S\in \mathcal{\tilde{C}}_K$.
\end{IEEEproof}
{As a remark, \cite[Thm~2]{Liu2011ControllabilityOC} associates the existence of a perfect matching in $\mathcal{H}_b(S)$ with a membership of $S$ to $\mathcal{C}_K$, whereas our result extends this previous result by associating the cardinality of a maximum matching in $\mathcal{H}_b(S)$ with a membership of $S$ to $\mathcal{\tilde{C}}_K$.} Since we invoke~\cite[Thm~2]{Liu2011ControllabilityOC}, our proposition also requires the accessibility condition. However, note that it is already satisfied by the strong connectivity of~$G$.
The proposition above also provides us with a systematic approach to calculate the smallest $K$ required for a nonempty $\mathcal{\tilde{C}}_K$, since we have $\emptyset \in \mathcal{\tilde{C}}_K$ if and only if $|\bar{m}(\emptyset)|\geq n-K$ holds. This method for finding the smallest $K$ coincides with the ones proposed in \cite[Thm~4]{pequito2015framework},~\cite[Thm~3]{pequito2013structured}.

\begin{example}
\label{exp: forwardFeas}
Consider a system described by $4$ nodes and the dynamic equations \eqref{eq: systemmodel} where
$${ A = \begin{bmatrix}
   0 & -0.5 & -0.8 & -0.6\\
   1 & 0 & 0 & 0\\
   1 & 0 & 0 & 0\\
   1 & 0 & 0 & 0
\end{bmatrix}.}$$
$G=(V,E)$ corresponding to this system is in Figure~\ref{fig:system}. For the metric in \eqref{eq:F(S)}, let $T=2$ and $\epsilon = 10^{-9}$. 

Consider the actuator placement on this system. We first study the minimum required cardinality for structural controllability. In the auxiliary bipartite graph $\mathcal{H}_b(\emptyset)$ shown in Figure \ref{fig:bipartite}, any maximum matching consists of 2 edges, that is, $|\bar{m}(\emptyset)|=2$. By Proposition \ref{prop: feasibility}, $\emptyset \in \mathcal{\tilde{C}}_K$ if and only if $\bar{m}(\emptyset)\geq 4-K+0$, that is, $K\geq2$. Therefore, we need at least 2 actuators to render the system structurally controllable.
Suppose $K=2$. The solution of the forward greedy is $\{v_3,v_4\}$. We depict the auxiliary bipartite graph $\mathcal{H}_b(\{v_3,v_4\})$ in Figure \ref{fig:bipartite} to check whether this actuator set is feasible. Maximum matching contains 4 edges, thus $\bar{m}(\{v_3,v_4\})=4\geq 4-2+2=n-K+k$. By Proposition~\ref{prop: feasibility}, $\{v_3,v_4\}$ belongs to $ \mathcal{\tilde{C}}_2.$

We now provide a counterexample based on the example above to show that the feasibility check method in \cite{clark2012leader} excludes feasible nodes from the consideration of the forward greedy algorithm.
Suppose $K = 3$. The feasibility check method in~\cite{clark2012leader} indicates that $\{v_1\}\notin \tilde{\mathcal{C}}_3$, because $v_1'$ is not missed by any maximum matching in $\mathcal{H}_b(\emptyset)$. However, since $\{v_3,v_4\}$ is structurally controllable, so is $\{v_1,v_3,v_4\}$. Then, $\{v_1\}\subset\{v_1,v_3,v_4\}$ implies that $\{v_1\}\in \tilde{\mathcal{C}}_3$.

\begin{figure}[t]
    \centering
    \includegraphics[width=0.35\linewidth]{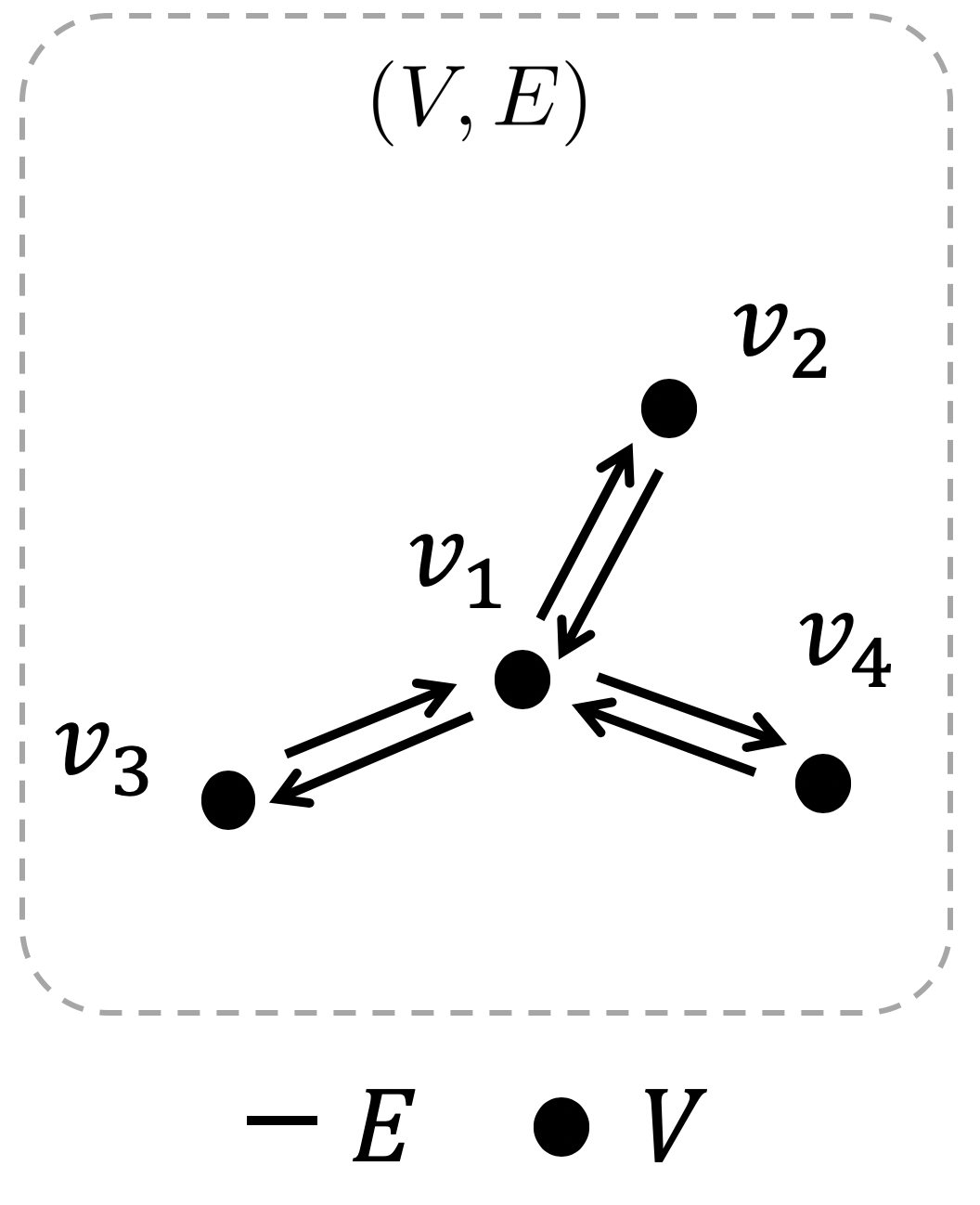}
    \caption{Graph for the 4-node system}
    \label{fig:system}
\end{figure}
\begin{figure}[t]
    \centering
    \includegraphics[width=0.68\linewidth]{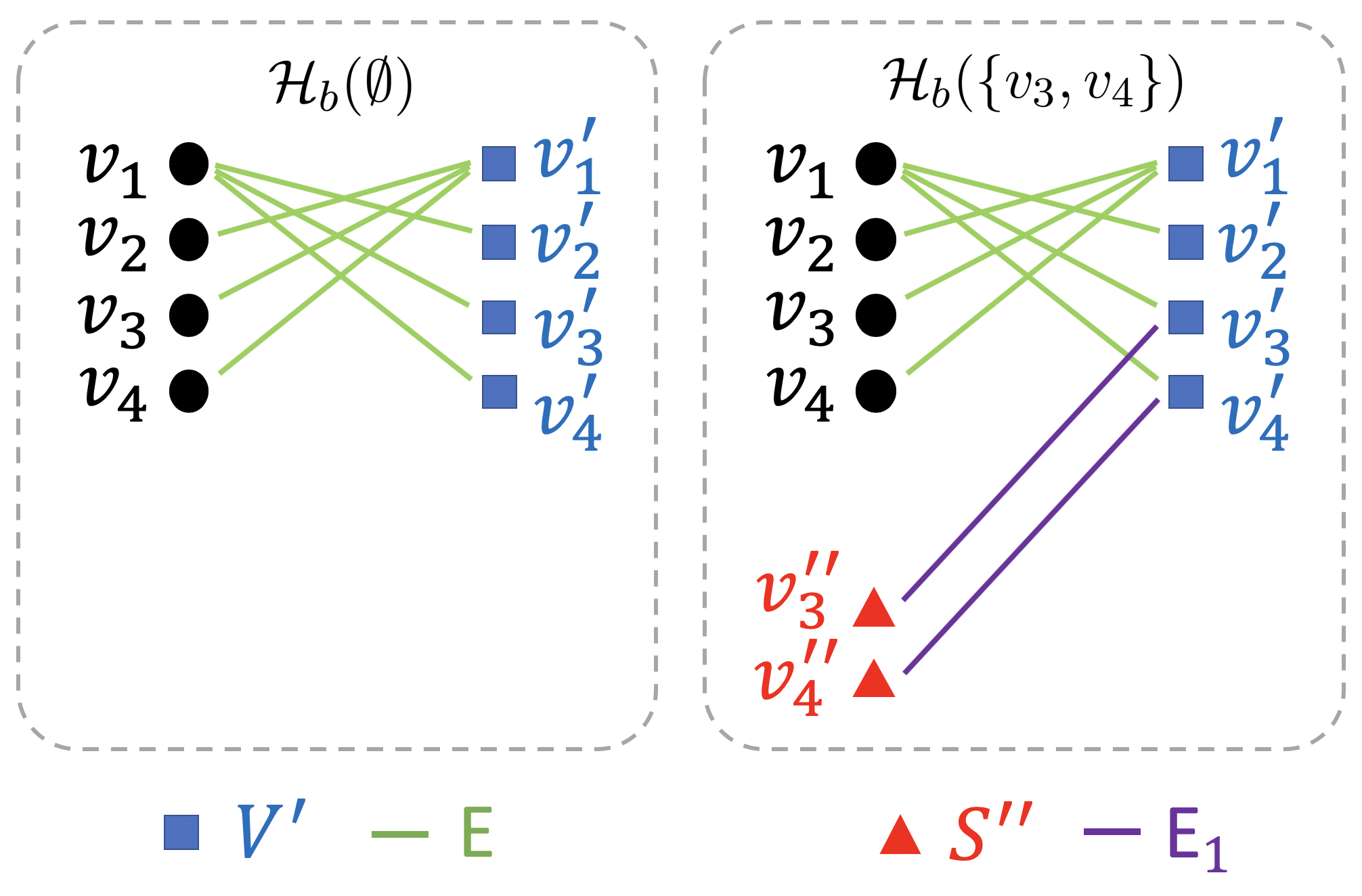}
    \caption{Auxiliary bipartite graphs $\mathcal{H}_b(\emptyset)$ and $\mathcal{H}_b(\{v_3,v_4\})$}
    \label{fig:bipartite}
\end{figure}
\end{example}

For our feasibility check, we still need a method to obtain a maximum matching in $\mathcal{H}_b(S)$. 
It is well-established that this can equivalently be done by solving a maximum flow problem~\cite{plummer1986matching}.
We refer to Appendix~\ref{apd: formulation of the max-flow problem} for details on formulating a maximum flow problem to obtain a maximum matching in $\mathcal{H}_b(S)$.
There are several algorithms for solving maximum flow problems. For instance, the Edmonds-Karp algorithm that we adopt in the numerical studies requires $O(pq^2)$ steps, where $p$ and $q$ respectively denote node cardinality and edge cardinality in the flow graph generated based on $\mathcal{H}_b(S)$\cite{Edmonds:1972:TIA:321694.321699}. For example, in $\mathcal{H}_b(\emptyset)$, $p = 2n+2$ and $q = 2n+|E|$. Thus, at each forward greedy iteration, we can examine in polynomial time whether $v\cup S^t$ belongs to $\mathcal{\tilde{C}}_K$ by finding the cardinality of the maximum matching in $\mathcal{H}_b(v\cup S^k)$.

\subsection{Feasibility check over $\mathcal{\tilde{R}}_K$}
\label{sec:reversefeasibilitycheck}

The reverse greedy algorithm has to determine whether $R\in \mathcal{\tilde{R}}_K$, or equivalently, whether any subset of the set $V\setminus R$ belongs to~${\mathcal{C}}_K$. Invoking the equivalence result of \cite[Thm~2]{Liu2011ControllabilityOC}, we can conclude that there exists a subset of $V\setminus R$ belonging to ${\mathcal{C}}_K$ if and only if there exists a perfect matching in $\mathcal{H}_b(V\setminus R)$ that covers at most $K$ elements of~${V''}\setminus {R''} $. This holds, since if every perfect matching in $\mathcal{H}_b(V\setminus R)$ covers $K+1$ or more nodes in $V''\setminus R''$, it would not be possible to find $K$ actuators from $V\setminus R$ satisfying structural controllability.\footnote{The feasibility check for the forward greedy does not limit the nodes in $S''$, because the cardinality of $S''$ is less than $K$ until the termination of the algorithm. Thus, the feasibility check of Section~\ref{sec:feasibilitycheck} is not applicable to $\tilde{\mathcal{R}}_K$.}
\looseness=-1

Recall that a maximum matching can be computed via the maximum flow algorithm. Analogous to the previous section, we need a feasibility check method for $\tilde{\mathcal{R}}_K$ by the means of the flow theory. We refer to Appendix~\ref{apd: formulation of the max-flow problem} for the preliminaries regarding flows in graphs. We first build an auxiliary graph, denoted by $\mathcal{H}_r(S)$, containing all the nodes in $\mathcal{H}_b(S)$. We let $s$ and $t$ be the sink and source of the flow, respectively. In addition, we add node $s''$ to $\mathcal{H}_r(S)$, which will enable encoding the cardinality limit on $V''\setminus R''$.
The edge set in $\mathcal{H}_r(S)$ is the union of three sets, $\mathsf{E}_b^f$, $\mathsf{E}^f_s$ and $\mathsf{E}^f_t$, all of which are directed. 
The edge set $\mathsf{E}^f_b$ is a copy of $\mathsf{E}\cup\mathsf{E}_1$, originally from $\mathcal{H}_b(S)$, but directed from $V\cup {S''}$ to $V'$ in $\mathcal{H}_r(S)$. The edge set $\mathsf{E}^f_s$ consists of edges from $s$ to all the nodes in $V$ and from~${s''}$ to all the nodes in~$S''$ along with edge of from~$s$ to~$s''$. Finally, the edge set $\mathsf{E}^f_t$ is composed of edges from all the nodes in $V'$ to $t$. 
All the edges have unit capacity except the edge from $s$ to $s''$ which has a capacity of $K$. Utilizing the graph $\mathcal{H}_r(S)$, we have the following proposition.
\begin{proposition}
\label{prop: feasibilityreverse}
Given the cardinality limit $K$ and an exclusion set $R$, we have $R\in \mathcal{\tilde{R}}_K$ if and only if there exists a flow~$g$ in $(\mathcal{H}_r({V}\setminus {R} ),c,s,t)$ with {\rm val}$(g)=n$.
\end{proposition}
\begin{IEEEproof}
From the definition of $\mathcal{\tilde{R}}_K$, $R\in \mathcal{\tilde{R}}_K$ is equivalent to the existence of $S\subset V\setminus R$ such that $S\in \mathcal{{C}}_K$. Via \cite[Thm~2]{Liu2011ControllabilityOC}, we know that these two conditions are equivalent to the existence a perfect matching in $\mathcal{H}_b(V\setminus R)$ that covers at most $K$ elements of $V''\setminus R''$. For the following, we prove that this equivalent condition holds if and only if there exists a flow~$g$ in $(\mathcal{H}_r({V}\setminus {R} ),c,s,t)$ with {\rm val}$(g)=n$.

``$\Rightarrow$": Given the perfect matching $m^*$, we use $\tilde{n}$ to denote the number of the elements in ${V''}$ that are adjacent to $m^*$. Clearly, $\tilde{n}\leq K$. For the following, we build the flow $g$ as a function of the edges in $\mathcal{H}_r({V}\setminus R)$. Suppose $m^f$ is a subset of $\mathsf{E}_b^f$ that corresponds with $m^*$ in $\mathsf{E}\cup \mathsf{E}_1$. We let $g(e)=1$ if $e$ belongs to $m^f$, $g((s,v))=1$ for any $v\in V$ incident with $m^f$, $g((s'',{v''}))=1$ for any ${v''}\in {V''}$ incident with $m^f$, $g(e_t)=1$ for any $e_t\in \mathsf{E}^f_t$ and finally $g((s,s''))=\tilde{n}$. It is easy to check $g$ is in fact a flow in $(\mathcal{H}_r({V}\setminus {R} ),c,s,t)$ with {\rm val}$(g)=n$.

``$\Leftarrow$": Let $\mathsf{E}^m = \{e^f\in \mathsf{E}^f_b \text{ }|\text{ }g(e^f)=1\}$. Then, define $m^*=\{e\in \mathsf{E}\cup \mathsf{E}_1\text{ }|\text{ } \exists e^f\in \mathsf{E}^m \text{ where $e^f$ is a copy}$ $\text{of $e$} \}$. It follows from {\rm val}$(g)=n$ that $m^*$ is a perfect matching in $\mathcal{H}_b(V\setminus R)$. Since the capacity limits are satisfied by the flow $g$, there are no more than $K$ elements in $V''\setminus R''$ covered by the perfect matching $m^*$.
\end{IEEEproof}


\begin{example}
We apply the reverse greedy algorithm to the system studied in Example \ref{exp: forwardFeas} with $K=2$. The first node excluded is $v_1$. To see that $v_1 \in \mathcal{\tilde{R}}_2$, we depict $\mathcal{H}_r(V\setminus v_1)$ in Figure \ref{fig:REVERSEgraph}. The maximum flow has a value of~4. Invoking Proposition \ref{prop: feasibilityreverse}, we conclude that $v_1$ belongs to $\mathcal{\tilde{R}}_2$.
\begin{figure}[t]
    \centering
    \includegraphics[width=0.58\linewidth]{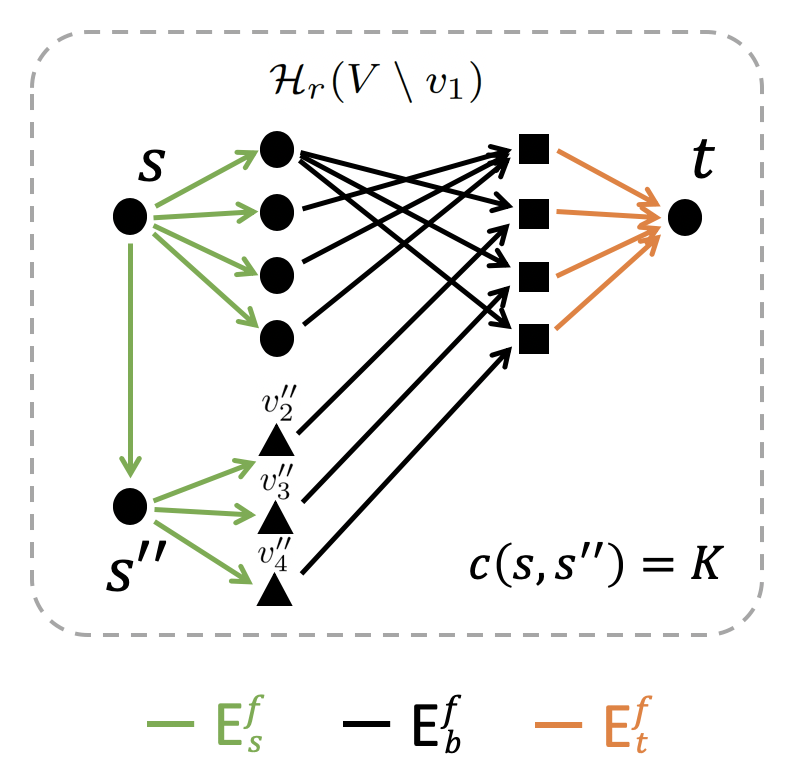}
	    \caption{Auxiliary graph $\mathcal{H}_r(V\setminus v_1)$}
    \label{fig:REVERSEgraph}
\end{figure}
\end{example}

Similar to Section \ref{sec:feasibilitycheck}, we adopt Edmonds-Karp algorithm to solve the maximum flow problem in the numerical studies. The algorithm requires $O(pq^2)$ steps where $p$ and $q$ are respectively the node cardinality and the edge cardinality of the flow graph $\mathcal{H}_r(V\setminus R)$. For example, in $\mathcal{H}_r(V)$, $p=3n+3$ and $q = 3n+|E|+1$, where $E$ is the edge set of $\mathcal{H}_r(V)$. Thus, at each reverse greedy iteration, we can examine in polynomial time whether $v\cup R^t$ belongs to $\mathcal{\tilde{R}}_K$.
\begin{remark}
{Greedy algorithms can also be applied when strong connectivity assumption on $G$ is relaxed.
Suppose $G$ is not strongly connected but it can be decomposed as $\cup^l_{i=1}G_{i}$, where $G_i$ is strongly connected for any $i$. 
In this case, if we suppose at least one actuator is already chosen and assigned for each subgraph~$G_i$, we would attain the accessibility condition discussed in Footnote~\ref{foot:9}. It is then possible to invoke \cite[Thm~2]{Liu2011ControllabilityOC} for the equivalence between perfect matching and structural controllability. This would make it possible to extend the proofs of Propositions~\ref{prop: feasibility} and~\ref{prop: feasibilityreverse} for the case when we are assigning additional actuators in such graphs. We kindly refer the readers to \cite{Liu2011ControllabilityOC,clark2017submodularity} for a detailed discussion on this condition. Similar to our paper, many works in the literature, such as~\cite{clark2012leader,clark2017submodularity}, assume that $G$ is strongly connected such that the accessibility condition is automatically attained.}
\looseness=-1
\end{remark}
\section{Numerical Results}
\label{sec:num}
In this section, we apply the greedy algorithms to problems based on randomly generated networks and a large power grid. 
All problems are solved on a computer equipped with 8 GB RAM and a 2.7 GHz dual-core Intel i5 processor. 

\subsection{Experiment on a 23-node network}
\label{sec:experiment}

\begin{figure}[t]
    \centering
    \includegraphics[width=0.88\linewidth]{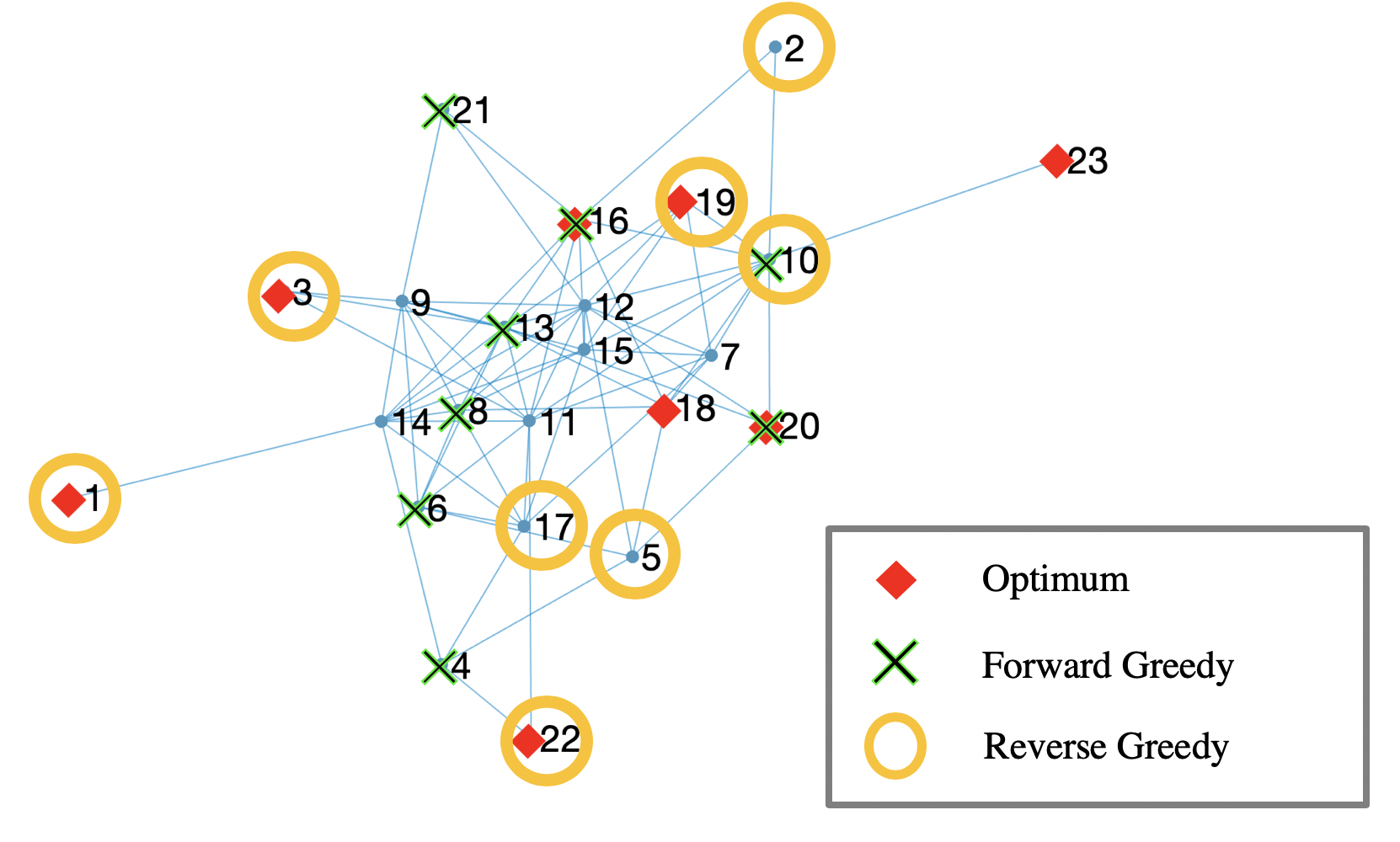}
    \caption{Greedy selection versus the optimal selection}
    \label{fig:experimentresult}
\end{figure}

We study a system model based on an undirected unweighted graph given in Figure~\ref{fig:experimentresult} generated via Octave Networks Toolbox~\cite{NetworkToolbox}. Different degrees are assigned to each vertex such that we can compare the sets $S^\mathsf{f}_\epsilon$ and $S^\mathsf{r}_\epsilon := V\setminus R^\mathsf{r}_\epsilon$ in terms of node connectivity. Specifically, vertex $i$ has a degree of $i$ if $i<12$ and a degree of $24-i$ if $i\geq12$. If there is an edge between vertex $i$ and $j$, we set $(A)_{ij}=(A)_{ji}=1$, otherwise the corresponding entries are $0$.\footnote{{This example involves a symmetric matrix, however our results are applicable for general matrices. In our future work, we aim to study the restricted concept of symmetric structural controllability~\cite{romero2018actuator}.}}

Let $T=1$ and $K=8$. We then apply Algorithm~\ref{alg:epsilon} to obtain a proper parameter $\epsilon$ for the forward greedy algorithm. We set $\xi=2$ and $\epsilon_0=10^{-3}$ arbitrarily. The actuator set returned in the first iteration is denoted by $S^\mathsf{f}_{\epsilon_0}$. The minimum eigenvalue corresponding to $W_T(S^\mathsf{f}_{\epsilon_0})$ is $\lambda_1(S^\mathsf{f}_{\epsilon_0}) = 1.9\times 10^{-4}$. Since $\epsilon_0>2\lambda_1(S^\mathsf{f}_{\epsilon_0})$, we continue with the second iteration. Let $\epsilon_1=\lambda_1(S^\mathsf{f}_{\epsilon_0})=1.9\times 10^{-4}$, we now have $S^\mathsf{f}_{\epsilon_1} = \{4,6,8,10,13,16,20,21\}$ and $\lambda_1(S^\mathsf{f}_{\epsilon_1}) = 2.0\times10^{-4}>\xi^{-1}\epsilon_1$. Thus, we can terminate the algorithm and pick $\epsilon_\mathsf{f} = \epsilon_1$ for the forward greedy algorithm. Using the same procedure, we obtain $\epsilon_\mathsf{r} = 1.4\times 10^{-4}$ for the reverse greedy algorithm. In this case, the solution is $S^\mathsf{r}_{\epsilon_\mathsf{r}} = V\setminus R^\mathsf{r}_{\epsilon_\mathsf{r}} = \{1,2,3,5,10,17,19,22\}$.

To assess the optimality of the sets $S^\mathsf{f}_{\epsilon_\mathsf{f}}$ and $S^\mathsf{r}_{\epsilon_\mathsf{r}}$, we generate the optimal solution $S^* = \{1,3,16,18,19,20,22,23\}$ by enumerating all feasible solutions. The average energy consumptions for all actuator sets are given by $F(S^\mathsf{f}_{\epsilon_\mathsf{f}}) = 9226.5$, $F(S^\mathsf{r}_{\epsilon_\mathsf{r}}) = 12126.2$, and $F(S^*) = 6052.7$. For this example, the forward greedy algorithm returns a better solution than the reverse greedy algorithm. Later, in randomized examples we see that this is not generally the case. 

Next, we analyze the performance guarantees in~\eqref{eq:estimationOFFepsilon} and \eqref{eq: performanceguaranteeforreverse} under the sets $S^\mathsf{f}_{\epsilon_\mathsf{f}}$ and $S^\mathsf{r}_{\epsilon_\mathsf{r}}$. For the forward greedy algorithm, we computed $F_{\epsilon_\mathsf{f}}(S^\mathsf{f}_{\epsilon_\mathsf{f}})= 6188.2 = 0.66 F(S^\mathsf{f}_{\epsilon_\mathsf{f}})$. The greedy submodularity ratio for the forward greedy algorithm is computed as $\gamma_{\epsilon_\mathsf{f}}^\mathsf{fg} = 1$. Then, we obtain
$F(S^\mathsf{f}_{\epsilon_\mathsf{f}}) = 1.5F_{\epsilon_\mathsf{f}}(S^\mathsf{f}_{\epsilon_\mathsf{f}})
        \leq 0.75F_{\epsilon_\mathsf{f}}(\emptyset)+0.75F_{\epsilon_\mathsf{f}}(S^*)=9.2\times10^{4}+0.75F_{\epsilon_\mathsf{f}}(S^*).$
In this example, the appearance of $F_{\epsilon_\mathsf{f}}(\emptyset)$ in the performance guarantee undermines its tightness. On the other hand, for the reverse greedy algorithm, the greedy submodularity ratio of the objective function $F_{\epsilon_\mathsf{r}}^\mathsf{r}$ is computed as $\gamma_{\epsilon_\mathsf{r}}^\mathsf{rg}<0.01$. This value is negligibly small making the performance guarantee in \eqref{eq: performanceguaranteeforreverse} loose.

\subsection{Node connectivity analysis on the sets $S^\mathsf{f}_{\epsilon_{\mathsf{f}}}$ and $S^\mathsf{r}_{\epsilon_{\mathsf{r}}}$} 

To gain additional insights into solution dependence on node connectivity, we now compare the greedy solutions with the optimal solutions in terms of the degrees of the selected actuators.
In the previous study, the forward greedy algorithm selects the actuator set $S^\mathsf{f}_{\epsilon_\mathsf{f}}$ in the order of 16, 13, 5, 8, 6, 20, 10, 21. In this sequence, the first four nodes feature high degrees. This is because the high degree nodes generally result in larger marginal gains at the earlier stages of the forward greedy algorithm. Let $d_\Sigma(S)$ denote the sum of the degrees of all the nodes in set $S$. Observe that $d_\Sigma(S^\mathsf{f}_{\epsilon_\mathsf{f}})=54$, whereas $d_\Sigma(S^\mathsf{r}_{\epsilon_\mathsf{r}})=35$ and $d_\Sigma(S^*)=30$. This demonstrates that the reverse greedy algorithm does not have a tendency to pick high degree nodes. We illustrate these sets in Figure \ref{fig:experimentresult}. 

To show that this observation is not restricted to this specific example, we build 20 random graphs with 23 nodes using Octave~\cite{NetworkToolbox}. These graphs are built as follows. For $i=1,4,7,$ node $i$, node $i+1$ and node $i+2$ have randomized degrees between $i$ and $i+2$. Node $10$ and node $11$ have randomized degrees between $10$ and $11$. Node $12$ has exactly 12 neighbors. For $i>12$, Node $i$ has a degree number the same as that of Node $24-i$. For each algorithm run, a proper parameter $\epsilon$ is picked via Algorithm \ref{alg:epsilon}.

Comparisons of different actuator sets can be found in Tables~\ref{tab: greedyalgorithmsonrandomgraphsdgreesun} and~\ref{tab: greedyalgorithmsonrandomgraphs}. The set $S_o^*$ refers to the best solution out of $1\times10^4$ random selections of cardinality $K=8$, while it is not computationally feasible to obtain the exact optimal solution for each case. 
Table~\ref{tab: greedyalgorithmsonrandomgraphsdgreesun} shows that the forward greedy generally yields an actuator set with a high degree sum when compared to the other solutions. {Finally, Table~\ref{tab: greedyalgorithmsonrandomgraphs} shows that in several cases the set returned by the forward greedy results in significantly worse value in the objective than the other two solutions. Generally, both greedy algorithms achieve comparable performance, and at the same time much better performance than what the theoretical guarantees suggest (as in Section~\ref{sec:experiment}). Thus, it could be useful to implement both polynomial-time algorithms, and choose the best solution.}

The total computation time for 20 forward greedy algorithm runs is 8205.0 seconds, whereas the time for 20 reverse greedy algorithm runs is 665.0 seconds. It turns out that for this problem the reverse greedy algorithm requires fewer queries to the computationally expensive feasibility check problem when compared to the forward greedy algorithm.   

\begin{table}[t]
	\setlength\tabcolsep{2pt} 
\centering
\caption{Degree sum comparisons of different solutions}
\resizebox{.48\textwidth}{!}{\begin{tabular}{|c|c|c|c||c|c|c|c|c|}
\hline
 Netw. & $d_\Sigma(S^\mathsf{f})$ & $d_\Sigma(V\setminus R^\mathsf{r})$ & $d_\Sigma(S^*_o)$  & Netw. & $d_\Sigma(S^\mathsf{f})$ & $d_\Sigma(V\setminus R^\mathsf{r})$ & $d_\Sigma(S^*_o)$    \\ \hline
 1 & 43 & 46 & 46
  & 11 & 41 & 38 & 42 \\ \hline
 2& 47 & 47 & 49
  & 12 & 60 & 44 & 53
   \\ \hline
 3& 48 & 34 & 41
  & 13 & 50 & 47 & 47
  \\ \hline
 4& 49 & 35 & 43
   & 14 & 55 & 49 & 41
    \\ \hline
 5 & 45 & 43 & 32
  & 15 & 53 & 41 & 43
    \\ \hline
 6& 58 & 39 & 41
   & 16 & 56 & 43 & 46
       \\ \hline
 7& 44 & 45 & 49
   & 17 & 48 & 52 & 70
    \\ \hline
 8 & 52 & 25 & 44
    & 18 & 57 & 48 & 54
      \\ \hline
 9& 55 & 37 & 45
   & 19 & 76 & 60 & 50
 \\ \hline
 10& 61 & 38 & 59
    & 20 & 57 & 43 & 44
     \\ \hline
 
\end{tabular}}
\label{tab: greedyalgorithmsonrandomgraphsdgreesun}
\end{table}

\begin{table}[t]
\setlength\tabcolsep{2pt} 
\centering
\caption{Metric comparisons of different solutions}
\resizebox{.48\textwidth}{!}{\begin{tabular}{|c|c|c|c||c|c|c|c|c|}
\hline
 Netw. & $F(S^\mathsf{f})$ & $F(V\setminus R^\mathsf{r})$ & $F(S^*_o)$  & Netw. & $F(S^\mathsf{f})$ & $F(V\setminus R^\mathsf{r})$ & $F(S^*_o)$  \\ \hline
 1& 66252 & 13493 & 11327  & 11 & 11981 &14993 &7300  \\ \hline
 2& 11527 &	10035 &	8461  & 12 & 9950 &	9388 &	9358   \\ \hline
 3& 15398 &	17371 &	8461  & 13 & 32212 & 6684 & 9364  \\ \hline
 4& 13679 & 15690 & 7593   & 14 & 9683 & 10681 & 6804    \\ \hline
 5& 14835 & 14430 & 8406  & 15 & 11235 & 8540 & 7188    \\ \hline
 6& 18207 & 15176 & 7870   & 16 & 14114 & 12004 & 6795 
       \\ \hline
 7& 8980 & 12650 & 10515   & 17 & 8658 &	9163 & 8416    \\ \hline
 8& 22600 & 26324 & 6638    & 18 & 8717 & 10587 & 6838      \\ \hline
 9& 10633 & 10483 & 9690   & 19 & 13760 & 9818 & 9336 \\ \hline
 10& 9676 & 13079 & 6173    & 20 & 10044 & 13169 & 9264     \\ \hline
 
\end{tabular}}
\label{tab: greedyalgorithmsonrandomgraphs}
\end{table}

\subsection{{Power electronic actuator placement for 118-bus system}} 

{We illustrate our results by placing power electronic actuators that can modulate power injections in the IEEE 118-bus test system, provided in~\cite{christie2000power}. Similar to \cite[\S 4.B]{summers2015submodularity}, each bus is assumed to follow the linearized swing equations, that is,
$M_i\ddot{\theta}_i+D_i\dot{\theta}_i = P_i-\sum_ja_{ij}(\theta_i-\theta_j)$ for all $i\in\{1,\ldots,118\}$,
where $P_i$ is the net power injection and $a_{ij}$ is characterized by the line parameters. We kindly refer to \cite{bergen1981structure} for the modeling details. If bus~$i$ is not in our actuator set, we have $P_i=0$. 
Buses without generators are assumed to have no inertia (even when they have loads connected), and they have a one-dimensional state, since $M_i=0$ (inertia). Buses with generators have inertias, and they are instead associated with a two-dimensional state vector that includes both $\theta_i$~and~$\dot{\theta}_i$. We highlight that each state corresponds to a new node in our system graph $G$, which in some sense represents an extended version of the original 118-bus power network.} 
\looseness=-1

{Our goal is to choose $K=50,\,70$ buses out of $118$ buses to inject power to minimize the average energy consumption given by \eqref{eq:F(S)} ($T=1$), while ensuring structural controllability. We take into account the fact that it is not possible to actuate some of the nodes of our system graph~$G$ by excluding them from the actuator sets and the feasibility check methods. These nodes correspond to ${\theta}_i$ originating from buses with inertia, since the dynamics $d\theta_i/dt=\dot{\theta}_i$ cannot be actuated.}

{Next, we implement the greedy algorithms. The results are shown in Table \ref{tab: 118-bus-test}, where $S^*_o$ is the best structurally controllable solution out of $10^5$ random selections with cardinality $K$. The parameters for both greedy algorithms can also be summarized as follows: $\epsilon = 1\times 10^{-10}$ for $K=50$ and $\epsilon = 1\times 10^{-8}$ for $K=70$, chosen arbitrarily. Both greedy algorithms perform significantly better than random selections. Note that there are $6.2\times10^{33}$ and $3.2\times10^{33}$ possible combinations to check for $K=50$ and $K=70$, respectively. Since the optimal solution is computationally out of reach, we will not analyse the performance guarantees in~\eqref{eq:estimationOFFepsilon} and \eqref{eq: performanceguaranteeforreverse} as we did in the previous sections.
In both cases $K=50$ and $K=70$, the reverse greedy performed slightly better than the forward greedy algorithm. Finally, as we expected, choosing a larger $K$ reduces the control cost for all three solution concepts. One can decide on $K$ by evaluating the overall cost reductions from the reductions in the metric and comparing them with the actuator installation costs.}
\looseness=-1

{\begin{table}[t]
\centering
\caption{Metric comparisons for the IEEE 118-bus test system}
\resizebox{.3\textwidth}{!}{\begin{tabular}{|c||c|c|c|}
\hline
 $K$ &   $F(S^{\mathsf{f}})$ &   $F(V\setminus R^{\mathsf{r}})$ & $F(S_{o}^*)$\\ \hline
 50&    $6.96\times10^7$ &  $3.01\times 10^7$ & $4.48\times 10^9$ \\ \hline
  70&  $1.94\times10^5$ &  $1.25\times 10^5$ & $4.43\times 10^6$ \\  \cline{1-4}
\end{tabular}}
\label{tab: 118-bus-test}
\end{table}
}

\section{Conclusions}

In this paper, our goal was to pick an actuator set to minimize a controllability metric based on average energy consumption while ensuring that the system is structurally controllable. To this end, we reformulated our problem as matroid optimization problems to apply both the forward and reverse greedy algorithms. For each algorithm, we provided a novel performance guarantee. For the implementation of the algorithms, we proposed feasibility check methods. In the numerics, we studied networks that are randomly generated based on degree lists. We observed that the forward greedy tended to select high-degree nodes in the early stages, whereas the overall performance of both algorithms were comparable.

Our future work involves exploiting the curvature of the objective function to derive a better performance guarantee for the forward greedy algorithm. We will exploit the problem structure to explain why algorithms performs significantly better than their performance guarantees. {We aim to investigate other structural controllability concepts from the literature.}

\bibliographystyle{IEEEtran}
\bibliography{IEEEabrv,library}

\appendix
\subsection{Performance guarantees from the literature}\label{app:tableg}
\bgroup
\def\arraystretch{1.5}
\begin{table}[t!]
\setlength\tabcolsep{2pt} 
\centering
\caption{Performance guarantees for the forward greedy algorithm applied to the maximization of increasing set functions}
\resizebox{.48\textwidth}{!}{\begin{tabular}{|c|c|c|c|c|}
\hline
 Setting & $\substack{\text{submodular}\\}$ & $\substack{\text{submodular}\\\text{with\ known\ curvature}}$ & $\substack{\text{known\ submodularity\ ratio}\\}$  & $\substack{\text{known\ submodularity\ ratio}\\\text{and\ known\ curvature}}$ \\ \hline
$\substack{\text{cardinality}\\\text{constraints}}$ & $1-e^{-1}$\cite{nemhauser1978analysis} & $(1-e^{-\alpha})/\alpha\text{\cite{conforti1984submodular}}$ & $1-e^{-\gamma}$\cite{bian2017guarantees}  & $(1-e^{-\alpha\gamma})/\alpha$\cite{bian2017guarantees} \\ \hline
$\substack{\text{general matroid}\\\text{constraints}}$ & 1/2\cite{fisher1978analysis} & $1/(1+\alpha)\text{\cite{conforti1984submodular}}$ & $\substack{\gamma^3/(1+\gamma^3)\ \text{in this paper}
\ \\ \text{see\ Appendix}~\text{\ref{apd: definitions of submodularity ratio}\ for\ \cite{pmlr-v80-chen18b}}}$  & $\substack{\text{open}\\}$  \\ \hline
\end{tabular}}
\label{tab:litrev}
\end{table}
\egroup 

{Before summarizing the relevant guarantees from the literature, we highlight that if we assign each actuator with a specific control cost, we would obtain a modular objective. There are many works that study modular objectives with or without structural controllability type constraints \cite{romero2018actuator,summers2015submodularity,doostmohammadian2018structural}. In case our objective is additive/modular, greedy algorithm on a matroid is known to always return an optimal solution~\cite{edmonds1971matroids}.}

Table~\ref{tab:litrev} summarizes the performance guarantees for the forward greedy algorithm applied to the maximization of increasing set functions. As a remark, \cite{bian2017guarantees} defines the curvature as we defined in Definition~\ref{def:curvature}, whereas \cite{conforti1984submodular} considers only 
the case where $S=\emptyset$ and $U=V\setminus \{v\}$ for all $v\in V$ and they call this notion the total curvature. In addition, the work of~\cite{sviridenko2017optimal} provides $1-(\alpha/e)$ guarantee that holds by high probability for maximizing increasing submodular functions with curvature over an arbitrary matroid constraint. This result relies on linear extensions of the objective by implementing the continuous greedy algorithm, which is subject to potential deviations from the guarantee due to the random rounding procedures. We refer to~\cite{sviridenko2017optimal} for other guarantees from the literature relying on variations of this algorithm. Finally, \cite{pmlr-v80-chen18b} relies on a randomized forward greedy providing guarantees only in expectation, see Appendix~\ref{apd: definitions of submodularity ratio}.
A very recent work~\cite{chamon2019matroid} provides $\gamma/2$ guarantee for the forward greedy algorithm applied to our setting with the submodularity ratio in Definition~\ref{def:submodularityratio}. This result does not rely on constructing a linear program. Instead, it relies on~\cite[Lem.~1]{chamon2019matroid}, which might have a mistake in its proof. Using the notation of \cite{chamon2019matroid}, let $\mathcal{G}_t$ be the greedy solution at iteration $t$, $\mathcal{X}^\star$ be the optimum, and $\mathcal{I}$ be the feasible region of the matroid. In the proof, the authors need to find an enumeration $\{x^\star_1,\ldots,x^\star_n\}=\mathcal{X}^\star$ such that $\mathcal{G}_t\cup \{x_{t+1}^\star\}\in\mathcal{I}$ for any $t$. Starting from $\mathcal{G}_0=\emptyset$, their inductive proof iteratively builds this enumeration by choosing $x_{t+1}^\star\in\mathcal{X}^\star$ with $\mathcal{G}_{t}\cup \{x_{t+1}^\star\}\in\mathcal{I}.$ The proof states ``$|\mathcal{G}_t|<|\mathcal{X}^\star|$ implies there exist an element $x_{t+1}^\star\in\mathcal{X}^\star$ such that $\mathcal{G}_t\cup \{x_{t+1}^\star\}\in\mathcal{I}$ (property~(iii) of Definition~\ref{def:matroid})". However, for their inductive proof to be correct, they instead have to prove that there exists an element $x_{t+1}^\star\in\mathcal{X}^\star\setminus\{x_{k}^\star\}_{k=1}^t$ such that $\mathcal{G}_t\cup \{x_{t+1}^\star\}\in\mathcal{I}$. We think that this statement could also be proved. In this case, their guarantee would be better. \looseness=-1

The guarantees above take as a reference objective evaluated at the empty set. The second part of this work studies the reverse greedy algorithm where the guarantees are with respect to the objective evaluated at the full set (see~\eqref{eq: performanceguaranteeforreverse} and~\eqref{eq: boundwithEMPTYSET}).

\subsection{Definitions of submodularity ratio}
\label{apd: definitions of submodularity ratio}

Let $\gamma_1$ denote the submodularity ratio of $f$ from Definition~\ref{def:submodularityratio}.
Observe that $\gamma = \gamma_1$ satisfies
\begin{equation}
\label{eq:definition for submodularity ratio2}
\gamma\rho_{U}(S) \leq \sum_{v\in U\setminus S}\rho_v(S),\ \forall S,U \subset V,
\end{equation}
{which can easily be obtained by decomposing the term on the left via telescoping sum.}
However, the largest $\gamma$ satisfying the above set of inequalities, denoted as $\gamma_2$, does not necessarily satisfy the inequalities in Definition~\ref{def:submodularityratio}. {This is true since the inequalities in \eqref{eq:definition for submodularity ratio2} can be regarded as a relaxation of~those in Definition~\ref{def:submodularityratio}.} Hence, we have $\gamma_2\geq\gamma_1$. There are previous studies in the literature defining the submodularity ratio as $\gamma_2$ instead of $\gamma_1$~\cite{bian2017guarantees,summers2017performance,pmlr-v80-chen18b}. In the proof of Theorem \ref{thm:upperlimitFormatroidconstraints}, as we are deriving (\ref{eq: submodularityratioTheorem}), we use the inequalities from Definition~\ref{def:submodularityratio}. One can verify that the inequalities in (\ref{eq:definition for submodularity ratio2}) would not allow us to derive (\ref{eq: submodularityratioTheorem}). Hence, the performance guarantee (\ref{eq: performanceguaranteeMatroidOpt}) does not extend to the submodularity ratio $\gamma_2$. 

In addition, the work of \cite{summers2017performance} obtains a lower bound for $\gamma_2$ for the metric $-F$ in \eqref{eq:F(S)} based on eigenvalue inequalities for sum and product of matrices. One can easily verify that this lower bound is also applicable to $\gamma_1$ from Definition~\ref{def:submodularityratio}.

\begin{figure}
    \centering
    \includegraphics[width=0.75\linewidth]{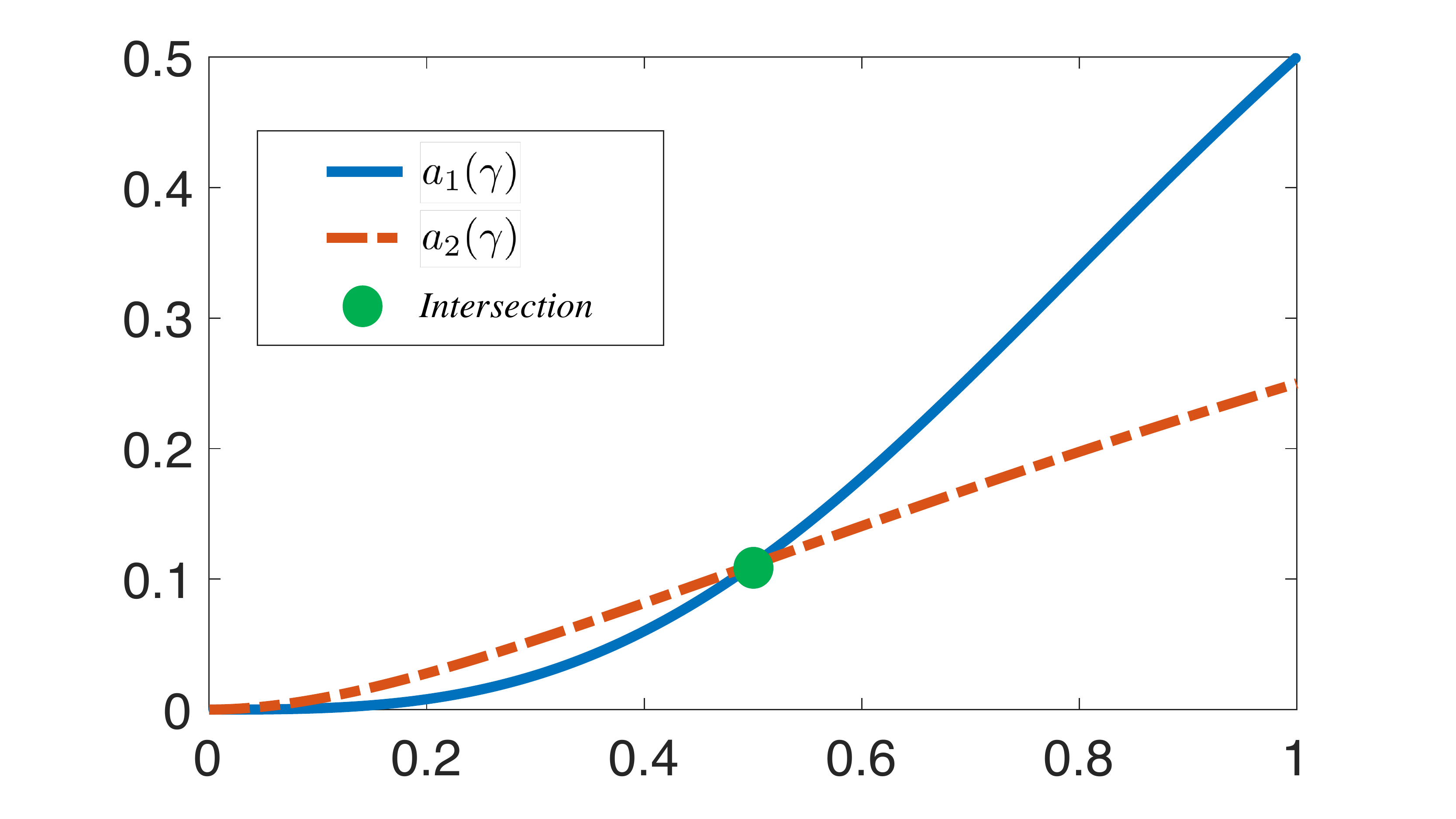}
    \caption{A comparison between two guarantees}
    \label{fig:comparisonbetweentwoguarantees}
\end{figure}

The work of \cite{pmlr-v80-chen18b} exploited the submodularity ratio defined by (\ref{eq:definition for submodularity ratio2}) and obtained a guarantee in expectation for the residual random (forward) greedy algorithm for matroid optimization problems featuring weakly submodular objective functions. We denote the final set returned by this algorithm as $S^{\text{RRG}}$. The guarantee provided in \cite{pmlr-v80-chen18b} for this class of randomized algorithms is $\frac{\mathbb{E}\left[f(S^{\text{RRG}})\right]-f(\emptyset)}{f(S^*)-f(\emptyset)}\geq \frac{\gamma_2^2}{(1+\gamma_2)^2}.$
Let $\gamma$ denote the theoretical lower bound derived in \cite{summers2017performance} for $-F$ in \eqref{eq:F(S)}. This lower bound satisfies $\gamma_2\geq\gamma_1\geq\gamma$. Since $\gamma$ is applicable to both \eqref{eq: performanceguaranteeMatroidOpt} and the guarantee in~\cite{pmlr-v80-chen18b}, {we let $a_1(\gamma) = \gamma^3/(1+\gamma^3)$ and $a_2(\gamma) = \gamma^2/(1+\gamma)^2$ denote the theoretical guarantees associated with Theorem~\ref{thm:upperlimitFormatroidconstraints} and the one in~\cite{pmlr-v80-chen18b}, respectively.} Two functions are plotted in Figure \ref{fig:comparisonbetweentwoguarantees}. The guarantee we derived in Theorem~\ref{thm:upperlimitFormatroidconstraints} is tighter than the one from~\cite{pmlr-v80-chen18b}, if the lower bound~$\gamma>0.5$ (it is also an ex-post guarantee).
\subsection{Proof of Lemma \ref{lmm: monotonicity}}
\label{apd: proofofmonotonicity}

(i) For any $S \subset V$ and any $v\in V\setminus S$, let $H(z)=(W_T(S)+zW_T(\{v\})+\epsilon I)^{-1}$. {Notice that
$\text{tr}(H(1)) = \text{tr}((W_T(S\cup \{v\})+\epsilon I)^{-1}) = F_{\epsilon}(S\cup \{v\}),$
since $W_T(S)=\int_{0}^{T} e^{A\tau}B(S) B^\top (S) e^{A^\top\tau} d\tau$ is additive, that is, $W_T(S)+W_T(\{v\})= W_T(S\cup \{v\})$.} Via the matrix inverse formula~\cite{matrixbook}, if $H(z)$ is invertible $\forall z\in(0,1)$, then $\text{tr}(H(z))$ is continuous and differentiable, and we have
$\frac{d(\text{tr}(H(z)))}{dz}=-\text{tr}(H(z)W_T(\{v\})H(z))< 0.$
This inequality holds since $H(z)$ is invertible and symmetric, and $W_T(\{v\})$ is positive semidefinite. Invoking the mean-value theorem, we have $\text{tr}(H(1))-\text{tr}(H(0))< 0$.

(ii) Recall from \eqref{eq:extentionofGammaAlphaDef} that the submodularity ratio of $-F_\epsilon$, denoted as $\gamma_\epsilon^\mathsf{f}$, satisfies
$
\gamma_\epsilon^\mathsf{f} = \min_{S,U,v\in V\setminus(S\cup U)}\frac{\rho_v(S)}{\rho_v(S\cup U)}. $
Since $-F_\epsilon$ is strictly increasing, ${\rho_v(S)}> 0$ and $\rho_v(S\cup U)> 0$ for any $v \in V \setminus (S\cup U)$. Thus $\gamma_\epsilon^\mathsf{f}>0$. \QEDA

\subsection{Proof of Proposition \ref{prop: matroid}}
\label{apd: matroid}
To prove this theorem, we show that given an actuator set $S$, structural controllability of $(A,B(S))$ can equivalently be formulated as structural controllability of the system with the set $S$ chosen as a leader set. Then, we use a result from~\cite{clark2012leader} showing the matroid structure of the structural controllability constraints in leader selection problems. This result builds on \cite{liu2011controllability}, which shows the equivalence between structural controllability and existence of a perfect matching in an auxiliary bipartite graph whenever the graph $G$ is strongly connected.
\looseness=-1

Define $N= V\setminus S$ and partition the state vector $x$ into $x_S$ and $x_N$. Our dynamics can equivalently be written as 
\begin{equation}
\label{eq: statepartition}
\begin{bmatrix}
    \dot{x}_{N} \\ \dot{x}_{S}
\end{bmatrix} = \begin{bmatrix}
 {A}_{NN} & {A}_{NS}  \\
 {A}_{SN} & {A}_{SS}
\end{bmatrix} \begin{bmatrix}
    x_{N} \\ x_{S}
\end{bmatrix} + \begin{bmatrix}
   0 & 0\\ 0 & I_{|S|}
\end{bmatrix}u,
\end{equation}
where $I_{|S|}\in \R^{|S|\times|S|}$ is the identity matrix.

In the leader selection problem, if the set $S$ is chosen as a leader set, it is assumed that the values of $x_S$ are directly dictated and are not influenced by the dynamics of $x_N$, see~\cite{patterson2010leader}.
Under this assumption, by treating $x_S$ as the input, the dynamics of $x_N$ are given by
$\dot{x}_N = A_{NN}x_N+A_{NS}x_S.$ Then, the leader set $S$ achieves structural controllability if  $(A_{NN},A_{NS})$ is structurally controllable, which would allow the values of $x_N$ to be steered to desired positions. Note that it is not clear whether we would achieve structural controllability when this is chosen as the set of actuators in actuator placement problem. 

From Definition \ref{def: structuralcontrollability}, the actuator set $S$ makes the system structurally controllable if and only if there exists a pair $(\hat{A},\hat{B})$ with the same structure as $(A,B(S))$ such that the controllability matrix $P\in \R^{n\times n^2}$,
$$P = \begin{bmatrix}
   0 & 0 & 0 &\hat{A}_{NS} & 0 & \hat{A}_{NN}\hat{A}_{NS}+\hat{A}_{NS}\hat{A}_{SS} &\cdots\\
   0 & I_{|S|} & 0 & \hat{A}_{SS} & 0 & \hat{A}_{SN}\hat{A}_{NS}+\hat{A}^2_{SS} & \cdots
   \end{bmatrix},$$
   has full rank.
Next, we claim that $P$ has full rank if and only if the following matrix $\tilde{P_1}\in \R^{|N|\times n^2}$ has full rank, 
$\tilde{P_1} = \begin{bmatrix}
   0 & 0 & 0 &\hat{A}_{NS} & \cdots &0 & \hat{A}^{j-1}_{NN}\hat{A}_{NS} & \cdots 
   \end{bmatrix}.$
To see this, notice that $P$ has full rank if and only if the submatrix $P_1\in \R^{|N|\times n^2}$ containing the first $|N|$ rows of $P$ has full rank. One can then show that there exists an upper triangular matrix $U\in \R^{n^2\times n^2}$ with unit diagonal entries such that $\tilde{P_1}= {P_1}U$.  Since $U$ is invertible, $\tilde{P_1}$ and $P_1$ have the same rank.
\looseness=-1

Then, we further claim that $\tilde{P_1}$ has full rank if and only if the following matrix $\bar{P_1}$ has full rank $\bar{P_1} = \begin{bmatrix}
\hat{A}_{NS} & \hat{A}_{NN}\hat{A}_{NS} & \cdots & \hat{A}^{|N|-1}_{NN}\hat{A}_{NS}\end{bmatrix}.$
Considering $|S|>0$ and thus $|N|-1\leq n-2$, for any $i>{|N|-1}$, $\hat{A}^{i}_{NN}\hat{A}_{NS}$ is in the span of the matrices $\hat{A}^{j}_{NN}\hat{A}_{NS}$, $j=\{0,1,\ldots,|N|-1\}$ by Cayley-Hamilton theorem. Hence, $\bar{P_1}$ has the same rank as $\tilde{P_1}$. This proves the claim.

In summary, $P$ has full rank if and only if $\bar{P_1}$ has full rank. By the definition of $\bar{P_1}$, $\bar{P_1}$ being full rank is equivalent to controllability of $(\hat{A}_{NN},\hat{A}_{NS})$. Hence, structural controllability of $(A,B(S))$ is equivalent to that of $({A}_{NN},{A}_{NS})$. 

Now, define $\mathcal{L}_K = \{S\text{ }|\text{ $|S| = K$ and } ({A}_{NN},{A}_{NS})$ is structurally controllable$\}$ and conclude that $\mathcal{L}_K = \mathcal{C}_K$. The set collection $\mathcal{L}_K$ consists of all the $K$ cardinality leader sets achieving structural controllability. From \cite[Thm 4]{clark2012leader}, we have that the pair $(V,\mathcal{\tilde{L}}_K)$, where $\mathcal{\tilde{L}}_K:=\{\Omega\text{ }|\text{ } \exists \text{ } S\in \mathcal{L}_K$ $ \text{such that }$ $ \Omega\subset S \}$, is a matroid if the graph $G$ is strongly connected.
Therefore, the pair $(V,\mathcal{\tilde{C}}_K)$ is also a matroid.
 \QEDA
\subsection{Proof of Theorem \ref{thm:upperlimitFormatroidconstraints}}
\label{apd: theorem1}
The idea of the proof extends the work in \cite{fisher1978analysis}, which derives a performance guarantee for matroid optimization featuring a submodular objective.
To assess the suboptimality of the actuator set $S^\mathsf{f}$, we need to find an upper bound for $f(S^*)-f(S^\mathsf{f})$. We denote $S^*=\{v^*_1,\ldots, v^*_K\}$ and notice
\begin{align}
f(S^*)-f(S^\mathsf{f}) &\leq f(S^*\cup S^\mathsf{f})-f(S^\mathsf{f})\nonumber\\
 &\hspace{-2.5cm}= \sum^{K}_{k=1} \rho_{v^*_k}(\{v^*_1,\ldots, v^*_{k-1}\}\cup S^\mathsf{f})\leq \gamma^{-1}\sum_{j\in S^*\setminus S^\mathsf{f}}\rho_j(S^\mathsf{f}),\label{eq: ideaTheorem}
\end{align}
where the first inequality is due to the monotonicity of~$f$ and the equality follows from a telescoping sum. The last inequality is from Definition \ref{def:submodularityratio}. To further bound $\sum_{j\in S^*\setminus S^\mathsf{f}}\rho_j(S^\mathsf{f})$, we have the following lemmas. For these lemmas, define $U^{-1}=\emptyset$, $U^K = V$, and
$s_{t} = |S^*\cap(U^{t+1}\setminus U^{t})|.$
\begin{lemma}
$\sum_{j\in S^*\setminus S^\mathsf{f}}\rho_j(S^\mathsf{f})\leq \gamma^{-1} \sum_{t=1}^{K}\rho_{t-1}s_{t-1}.$
\label{lmm: totalInequality}
\end{lemma}
\begin{IEEEproof}
From Definition \ref{def:submodularityratio}, we have
\begin{equation}
\begin{aligned}
\rho_j(S^\mathsf{f})  & \leq \gamma^{-1}\rho_j(S^{t-1}),\, \forall t\leq K\text{, } \forall j\in V.
\end{aligned}
\label{eq: submodularityratioTheorem}
\end{equation}
Since $U^{t_1}\subset U^{t_2}$ for any $t_1<t_2$, notice that
$V = U^K =  \bigcup^{K}_{t=0}(U^t\setminus U^{t-1}).$
Considering $U^{t_1}\setminus U^{t_1-1}$ and $U^{t_2}\setminus U^{t_2-1}$ are disjoint, we know that these sets constitute a partition of $V$.
Since there is no subset of $U^0$ belonging to $\mathcal{F}$, we have $S^*\cap U^0 = \emptyset$. Using the partition of $V$, we can partition $S^*$ as:
$S^*={{{\bigcup}^{{K}}_{t=1}}}(S^*\cap(U^t\setminus U^{t-1})).$
Combining this with (\ref{eq: submodularityratioTheorem}), 
{\medmuskip=.45mu\thinmuskip=.45mu\thickmuskip=.45mu\begin{equation}
\sum_{j\in S^*\setminus S^\mathsf{f}}\rho_j(S^\mathsf{f})  \leq \sum_{j\in S^*}\rho_j(S^\mathsf{f})= \sum_{t=1}^{K}\sum_{j\in S^*\cap(U^{t}\setminus U^{t-1})}\frac{1}{\gamma}\rho_j(S^{t-1}).
\label{eq:decomposition of S^*_m}
\end{equation}}Notice that all the nodes in $U^{t-1}$ have been considered by the feasibility check before $v^\mathsf{f}_t$. Since the greedy algorithm first checks the elements in $V\setminus U^{t-1}$ with larger marginal gains when added to $U^{t-1}$, we have that 
$\rho_{t-1} = \max_{j\in V\setminus U^{t-1}}\rho_j(S^{t-1}).$
Considering $V\setminus U^{t-1}= \cup_{i=t}^{K} (U^i\setminus U^{i-1})$,
for any $t'\geq t$,
\begin{equation}
\begin{aligned}
\rho_{t-1}\geq \rho_j(S^{t-1}), \forall j\in U^{t'}\setminus U^{t'-1}.
\end{aligned}
\label{eq: greedypickproperty}
\end{equation}
Thus, for any $j\in S^*\cap(U^t\setminus U^{t-1})$, we have $\rho_j(S^{t-1})\leq \rho_{t-1}$ and 
\begin{equation}
    \sum_{j\in S^*\cap(U^{t}\setminus U^{t-1})}\rho_j(S^{t-1}) \leq \rho_{t-1} s_{t-1}.
    \label{eq: backwardnumberinequality}
\end{equation}
Now combining (\ref{eq:decomposition of S^*_m}) and (\ref{eq: backwardnumberinequality}), it is straightforward that
$\sum_{j\in S^*\setminus S^\mathsf{f}}\rho_j(S^\mathsf{f})\leq  \sum_{t=1}^{K}\gamma^{-1}\rho_{t-1}s_{t-1}.$  
\end{IEEEproof}

\begin{lemma}\label{lmm: secondInequality}
For any $t\in \{1,\ldots,K\}$, we have $\sum_{i=1}^{t}s_{i-1} \leq t.$\end{lemma}

The above lemma is proven by \cite{fisher1978analysis} for $\gamma=1$, and it holds also when $\gamma\neq 1$ since its proof exploits only the matroid structure. 
The proof is included below for the sake of completeness.

\begin{IEEEproof}
We claim that any independent subset of $U^t$ has a cardinality at most $t$. Otherwise, due to $\mathcal{F}$ being a matroid, there exists $j\in U^t \setminus S^t$ $\text{such that } S^t\cup \{j\}$ is independent. Since $j\in U^t$ and $U^t = \cup^t_{i=0} (U^i\setminus U^{i-1})$ is a partition, there exists $t'\leq t$ such that $j\in U^{t'}\setminus U^{t'-1}$. Since $S^{t'}\cup\{j\}\subset S^{t}\cup\{j\}$, $S^{t'}\cup\{j\}$ is independent. By the mechanism of the greedy algorithm, we know $j$ passes the feasibility check ahead of $v^\mathsf{f}_{t'+1}$, which contradicts the fact that $j$ is discarded. Then, notice that $S^*\cap U^t$ is an independent subset of $U^t$. Hence, its cardinality is no more than $t$ according to the above claim. The partition $U^t = \cup^t_{i=0} (U^i\setminus U^{i-1})$ gives us that $\sum^t_{i=1}s_{i-1} = |S^*\cap U^t|\leq t.$
\end{IEEEproof}

We use Lemma~\ref{lmm: secondInequality} to obtain an upper bound to the right-hand side of Lemma~\ref{lmm: totalInequality} and consequently to derive an upper bound of $f(S^*)-f(S^\mathsf{f})$. The following explains these steps.

\begin{IEEEproof}[Proof of Theorem~\ref{thm:upperlimitFormatroidconstraints}] First, we consider the case in which $\rho_i$, $i = 0,\ldots, K-1$, are distinct. We define $t_1$ such that $\rho_{t_1-1}$ is the largest among $\rho_0, \rho_1,\ldots,\rho_{K-1}$ and $t_2$ such that $\rho_{t_2-1}$ is the largest among $\rho_{t_1}, \rho_{t_1+1},\ldots,\rho_{K-1}$. Following the same pattern we have $t_1,t_2,\ldots,t_p,$ where $t_p = K$. Since $s_{i}\geq 0$ is bounded by Lemma~\ref{lmm: secondInequality}, to give an upper bound to the right-hand side of Lemma~\ref{lmm: totalInequality}, we construct a linear program as follows,
{\medmuskip=.78mu\thinmuskip=.78mu\thickmuskip=.78mu\begin{equation}
 \max_{s_0,\ldots,s_{K-1}\geq 0}
 \sum_{i=1}^{K} \rho_{i-1}\,s_{i-1}\ \mathrm{s.t.}\, \sum_{i=1}^{t}s_{i-1} \leq t,\ t = 1,\ldots,K.
\label{eq:matroidconstrained optimization transferred}
\end{equation}}

Let $s_{i-1}^*$, $i = 1,2,\ldots,K,$ denote the optimal solution. We claim $s_{t_1-1}^* = t_1$. Otherwise, $s_{t_1-1}^* < t_1$ and due to Lemma~\ref{lmm: secondInequality} two situations might happen, \text{a)} $\sum^{t_1}_{i=1} s^*_{i-1} = t_1$ or \text{b)} $\sum^{t_1}_{i=1} s^*_{i-1} < t_1$.

For case \text{a)}, we obtain $\sum^{t_1-1}_{i=1} s^*_{i-1} > 0$. It follows that there exists $l < t_1\text{ such that } s^*_{l-1}>0$. Then, we decrease $s^*_{l-1}$ by $\delta>0$ and increase $s^*_{t_1-1}$ also by $\delta$. The value of $\delta$ is small enough so that $s^*_{l-1}>0$. This operation decreases $\sum^{t}_{i=1} s^*_{i-1}$ for $l\leq t\leq t_1-1$ and keeps the sum unchanged for any other $t$, so the constraints of (\ref{eq:matroidconstrained optimization transferred}) are not violated. Also considering that $\rho^*_{t_1-1}>\rho^*_{l-1}$, after these changes, the objective function is strictly greater than the value obtained at the original optimum. Thus, case \text{a)} is impossible. 
For case \text{b)}, we collect all the integers $l>t_1$ satisfying $s^*_{l-1}>0$. Assume $l_q>\cdots>l_1>t_1$. We have $q\geq 1$. Otherwise, $s^*_{l-1}=0$ for any $l>t_1$ and we can increase $s^*_{t_1-1}$ by a small amount to obtain a greater value of the objective without violating constraints. Knowing that $s^*_{l_1-1}>0$ and following the same reasoning provided for the case \text{a)}, we increase $s^*_{t_1-1}$ and decrease $s^*_{l_1-1}$ with the same amount. This way, an objective value is obtained larger than that evaluated at the original optimum. Thus, case \text{b)} is impossible. 

In conclusion, $s_{t_1-1}^* = t_1$ and (\ref{eq:matroidconstrained optimization transferred}) is equivalent to 
{\begin{equation}
\begin{aligned}
& \max_{s_{t_1},\ldots,s_{K-1}\geq 0}
& & \sum_{i=t_1+1}^{K} \rho_{i-1} s_{i-1}
\\ 
&\ \ \ \ \ \mathrm{s.t.} & & \hspace{-.4cm}\sum_{i=t_1+1}^{t}s_{i-1} \leq t-t_{1} , \text{ } t = t_1+1,\ldots,K.
\end{aligned}
\label{eq:matroidconstrained optimization transferred2}
\end{equation}}
We determine $s^*_{t_2-1}$ in the same way as we determine $s^*_{t_1-1}$ in (\ref{eq:matroidconstrained optimization transferred}). By repeating the above procedure we obtain the solution
  {\begin{equation}
    s_{i-1}^*=
    \begin{cases}
      t_1, & \text{if } i=t_1,\\
      t_{j}-t_{j-1}, & \text{if } i=t_j\text{ and $j\neq 1$},\\
    0, & \text{otherwise}.
    \end{cases}
    \label{eq:soptimal}
  \end{equation}}
  If $\rho_{i}$, $i=0,\ldots,K-1$ are not distinct and there exist $i_1<i_2<\cdots<i_q$ with $\rho_{i_1} = \rho_{i_2} = \cdots = \rho_{i_q}$. We let $s^*_{i_1} = \cdots =s^*_{i_{q-1}} = 0 $ and obtain the same solution as (\ref{eq:soptimal}). Next, notice
{\medmuskip=.4mu\thinmuskip=.4mu\thickmuskip=.4mu\begin{equation}
    \rho_{i_2} = f(S^{i_2+1})-f(S^{i_2})\leq \gamma^{-1}(f(S^{i_1}\cup v^\mathsf{f}_{i_2+1})-f(S^{i_1}))\leq \frac{\rho_{i_1}}{\gamma},
    \label{eq:rhoi1rhoi2}
  \end{equation}}where the first inequality comes from the definition of submodularity ratio, while the second is due to (\ref{eq: greedypickproperty}).
Substituting the optimal solution into the objective, considering (\ref{eq:rhoi1rhoi2}), we~have
 {\medmuskip=.78mu\thinmuskip=.78mu\thickmuskip=.78mu\begin{equation}
 \begin{aligned}
    \sum_{i=1}^{K} \rho_{i-1} s^*_{i-1} &= t_1\rho_{t_1-1}+\cdots+(t_p-t_{p-1})\rho_{t_{p}-1}\\
    &\hspace{-2cm}\leq \gamma^{-1}\sum^p_{k=1}\sum_{i=t_{k-1}+1}^{t_k}\rho_{i-1}=\gamma^{-1}\sum^{K}_{i = 1}\rho_{i-1}= \gamma^{-1}(f(S^\mathsf{f})-f(\emptyset)).
    \end{aligned}
    \label{eq:backwardInequalities}
  \end{equation}}
Combining (\ref{eq: ideaTheorem}), Lemma~\ref{lmm: totalInequality} and (\ref{eq:backwardInequalities}), we have
$f(S^*)-f(S^\mathsf{f}) \leq \gamma^{-1}\sum_{j\in S^*\setminus S^\mathsf{f}}\rho_j(S^\mathsf{f})\leq \gamma^{-2} \sum_{i=1}^{K} \rho_{i-1}s_{i-1}^* \leq \gamma^{-3} (f(S^\mathsf{f})-f(\emptyset)).$
By rewriting this, we have 
$\frac{f(S^\mathsf{f})-f(\emptyset)}{f(S^*)-f(\emptyset)}\geq\frac{\gamma^3}{\gamma^3+1}.$
\end{IEEEproof}
\subsection{Proofs of Propositions~\ref{prop:rechar},~\ref{prop: impossibilityResults1}, and~~\ref{prop: impossibilityResults2}}\label{app:props}
\begin{IEEEproof}[Proof of Proposition~\ref{prop:rechar}]
 Let $\rho_i(S) = -F_\epsilon(S\cup \{i\})-\left(-F_\epsilon(S)\right)$. Given $S$ and $U$, we denote $Y= U\setminus S$ and $R = V\setminus (S\cup U).$ Notice that if $i\notin S\cup U,$
 $\frac{\rho_i(S)}{\rho_i(S\cup U)} = \frac{F^\mathsf{r}_\epsilon(R\cup Y)-F^\mathsf{r}_\epsilon(R\cup Y\setminus\{i\})}{F^\mathsf{r}_\epsilon(R)-F^\mathsf{r}_\epsilon(R\setminus \{i\})}.$
 From Definitions \ref{def:submodularityratio} and \ref{def:curvature} we know that for all possible combinations of $S$ and $U$, the left-hand side has the least upper bound ${1}/({1-\alpha_\epsilon^\mathsf{f}})$ and the greatest lower bound $\gamma_\epsilon^\mathsf{f}$ while the right-hand side has the least upper bound ${1}/{\gamma_\epsilon^\mathsf{r}}$ and the greatest lower bound $1-\alpha_\epsilon^\mathsf{r}$. Consequently, we obtain $\gamma_\epsilon^\mathsf{r} = 1-\alpha_\epsilon^\mathsf{f}$ and $\alpha_\epsilon^\mathsf{r} = 1-\gamma_\epsilon^\mathsf{f}.$
\end{IEEEproof}
\begin{IEEEproof}[Proof of Proposition~\ref{prop: impossibilityResults1}]
We prove this by providing a counter-example where the greedy algorithm returns an arbitrarily poor solution as $\gamma$ goes to zero. We study a special instance of Problem~\eqref{eq:generalized reverse version} shown as
\begin{equation}
\underset{R\subset V}{\min}\ f(R)
\ \mathrm{s.t.}\
|R|\leq N\text{ and $|R|\geq N$,}
\label{eq:exampleNecessityofGamma}
\end{equation}
where $N=2$. Function $f$ is defined on the ground set $V = \{v_1,v_2,v_3\}$, with $f(\emptyset)=0$, $f(\{v_1\}) = \delta>0$, $f(\{v_2\}) = f(\{v_3\})= 2\delta$, $f(\{v_2,v_3\}) = 4\delta$, $f({\{v_1,v_2\}})=f(\{v_1,v_3\})=1$ and $f({\{v_1,v_2,v_3\}})$ $=2$. It is easy to verify that $f$ is strictly increasing, $\gamma = \delta/(2-4\delta)$, $\alpha = 0$ (supermodular) and $|R|\leq 2$ is a matroid constraint. The forward greedy algorithm would select $R^\mathsf{r} =\{v_1,v_2\}$ instead of $R^* = \{v_2,v_3\}$, thus 
$\frac{f(R^\mathsf{r})-f(\emptyset)}{f(R^*)-f(\emptyset)}  = \frac{1}{4\delta}.$
The value of $\delta$ can be chosen arbitrarily small to ensure $\gamma = \delta/(2-4\delta)$ goes to zero. This would imply that the set returned can be arbitrarily poor compared to the optimal solution in case we do not have a lower bound on the submodularity ratio.
\end{IEEEproof}
\begin{IEEEproof}[Proof of Proposition~\ref{prop: impossibilityResults2}]
We prove this by providing a counter-example for Problem \eqref{eq:exampleNecessityofGamma}. Let $\delta$ be a small number. Define $V=\{v_1,\ldots,v_n\}$ where $n>2N$. Let $\bar S=\{v_1,\ldots,v_N\}$. The function $f$ is given by 
$f(S)=\min\{1+|S\cap\bar S|,|S|\}+\delta|S|.$
This function is strictly increasing. Next, we show that the first term in~$f$, $g(S)=\min\{1+|S\cap\bar S|,|S|\}$ for all $S$, is submodular. Observe that the marginal gains are either $1$ or $0$. Whenever the marginal gain of adding an element to set $S$ is~$0$, the term on the left would be the minimum. This implies that if we add the same element to a superset of $S$, we would again obtain the marginal gain $0$. This concludes the submodularity of $g$. The function~$f$ is submodular, $\gamma =1$, since it is given by the sum of a modular and a submodular function. Finally, it can be verified that the function~$f$ has the curvature $\alpha  = 1/(1+\delta)$. 
A greedy solution is $R^\mathsf{r}=\{v_1,\ldots,v_N\}$, whereas an optimal solution is $R^*=\{v_n,\ldots,v_{n-N+1}\}$. Then, we have
$$\frac{f(R^\mathsf{r})-f(\emptyset)}{f(R^*)-f(\emptyset)} =\frac{f(\{v_1,\ldots,v_N\})}{f(\{v_n,\ldots,v_{n-N+1}\})} = \frac{N+\delta N}{1+\delta N}.$$
The value of $\delta$ can be chosen arbitrarily small to ensure $\alpha$ goes to $1$. This would imply that without an upper bound less than 1 on the curvature $\alpha$, we cannot obtain a performance guarantee better than $\lim_{\delta\rightarrow 0^+}\frac{N+\delta N}{1+\delta N}=N$.
\end{IEEEproof}

\subsection{Formulation of the maximum flow problem}
\label{apd: formulation of the max-flow problem}
Given a directed graph $D$ with two distinguished nodes $s$ (the \text{source}) and $t$ (the \text{sink}), denote the node set in this graph $D$ as $V(D)$ and denote the edge set in this graph $D$ as $E(D)$. Suppose no edge is directed into $s$ or out of $t$. Let $c:E(D)\to \mathbb{R}_+$ be a function that assigns to any edge $(u,v)$ in $E(D)$ a nonnegative value $c(u,v)$ called the \text{capacity} of the edge. Any function $g:E(D)\to \mathbb{R}_+$ is called a $\text{flow}$ in $(D,c,s,t)$ if it satisfies the following two conditions: \text{a)} $\sum_u g(u,w) = \sum_v g(w,v)$ for any $w\in V(D)\setminus \{s,t\}$ and $\text{b)}$ $g(u,v)\leq c(u,v)$ for all $(u,v)\in E(D)$. The first condition is the nodal balance, whereas the second condition is the capacity limits. The sum $\sum_wg(s,w)$ is called the \text{value} of the flow $g$ and denoted as $\text{val}(g)$. The maximum flow problem is formulated as finding the flow in $(D,c,s,t)$ with the maximum value. For the undirected bipartite graph $\mathcal{{H}}_b(S)$, we direct all the edges from $V$ to $V'$ and place two other nodes $s$ and $t$. Directed edges are built from $s$ to all the nodes in $V\cup V''$ and from all the nodes in $V'$ to $t$. Based on this new digraph $D$, we construct a capacity function such that $c(u,v)=1$ for any $(u,v)\in E(D)$. It is easy to verify that the maximum value of a flow in $(D,c,s,t)$ is equivalent to the cardinality of the maximum matching in $\mathcal{{H}}_b(S)$~\cite{plummer1986matching}.

\end{document}